\documentclass[a4paper]{article}

\usepackage[english]{babel}
\usepackage[utf8x]{inputenc}
\usepackage[T1]{fontenc}
\usepackage{lipsum}
\usepackage{blindtext}
\usepackage{graphicx,amssymb,amsmath,textcomp}
\usepackage{newpxtext} % ``New PX" font with original code \usepackage{newpxtext,newpxmath}

\usepackage[a4paper,top=3cm,bottom=2cm,left=3cm,right=3cm,marginparwidth=2cm]{geometry}
\usepackage{tikz,tikz-cd}
\usetikzlibrary{matrix,arrows,positioning,calc}
\usepackage{verbatim}

\usepackage[nottoc]{tocbibind} %Includes "References" in the table of contents
\usepackage{amsmath,amsthm}
\usepackage{faktor}
\usepackage{xfrac} 
\usepackage{graphicx}
\usepackage[colorinlistoftodos]{todonotes}
\usepackage[colorlinks=true, allcolors=blue]{hyperref}
\usepackage{mathtools}

\newcommand{\Hom}{{\rm Hom}}
\newcommand{\vdim}{{\rm vdim \ }}
\newcommand{\Ima}{{\rm Im}}

\newcommand{\rank}{{\rm rank \ }}
\newcommand{\cspec}{{\rm spec}}
\newcommand{\spec}{{\rm Spec}}
\newcommand{\Span}{{\rm Span}}

\newcommand{\N}{{\mathbb N}}

\newcommand{\K}{{\mathbb K}}

\newcommand{\dR}{{d_{dR}}}
\newcommand{\bfem}[1]{\textbf{\emph{#1}}}

\title{\textbf{On Shifted Contact Derived Artin Stacks}}
\author{Kadri İlker Berktav \footnote{Bilkent University, Department of Mathematics, Ankara, Turkey;  e-mail: kadri.berktav@bilkent.edu.tr.
	} 
} 

\date{\vspace{-5ex}}

\begin{document}
	
%%%%%%%%%%%%%%%%%%%% Text italic %%%%%%%%%%%%%%%%%%%%%%%%%%%%
\theoremstyle{plain}
\newtheorem{theorem}{Theorem}[section] %[section] or [subsection]
\newtheorem{lemma}[theorem]{Lemma} 
%NOT: araya [theorem] yazdigimiz zaman Lemma/propositon etc.'lari kendi icinde siralamak yerine 2. kod ile verilen genel numaralandirmayi (section veya subsection) takip eder.
\newtheorem{proposition}[theorem]{Proposition} %[section] or [subsection]
\newtheorem{corollary}[theorem]{Corollary}%[section]

%%%%%%%%%%%%%%%%%%%% Text roman %%%%%%%%%%%%%%%%%%%%%%%%%%%%%
\theoremstyle{definition}
\newtheorem{notations}[theorem]{Notations}
\newtheorem{notation}[theorem]{Notation}
\newtheorem{remark}[theorem]{Remark}
\newtheorem{observation}[theorem]{Observation}
\newtheorem{definition}[theorem]{Definition}
\newtheorem{condition}[theorem]{Condition}
\newtheorem{construction}[theorem]{Construction}
\newtheorem{example}[theorem]{Example}%[section]

\let\pf\proof
\let\epf\endproof
\numberwithin{equation}{section} %or {subsection}

\maketitle

\begin{abstract}
 This is a sequel of  \cite{kib} % our previous work arXiv:2209.09686 
  on the development of derived contact geometry. In \cite{kib}%that work
  , we formally introduced {shifted contact structures} on derived stacks. We then gave a {Darboux-type theorem} and the notion of {symplectification}  \emph{only for}  $k$-shifted contact derived schemes, with $k<0$.
 
 In this paper, we  extend the results of \cite{kib} %the work arXiv:2209.09686
 from derived schemes to derived Artin stacks and provide some examples of contact derived Artin stacks. In brief, we first show that for $k<0$, every $k$-shifted contact derived Artin stack admits a contact Darboux atlas. Secondly, we canonically describe  the symplectification of a  derived Artin stack equipped with a $k$-shifted contact structure, where $k<0$. Lastly, we give several constructions of contact derived stacks using certain cotangent stacks and shifted prequantization structures.
 \end{abstract}

\tableofcontents

    \section{Introduction and summary}

   As a relatively popular mainstream in the literature, generalized versions of certain familiar geometric structures have been introduced and studied in the context of derived algebraic geometry. 
   For example,  \cite{PTVV,CPTVV} focus on $k$-shifted Symplectic and Poisson geometries. Furthermore, \cite{Brav,JS,BenBassat} provide some applications and local constructions. 
   
   In \cite{kib}, on the other hand, we formally described shifted contact structures on derived (Artin) stacks and investigated some interesting consequences, such as  a Darboux-type theorem and the  process of canonical symplectification for negatively shifted contact derived schemes. These observations essentially motivate us to make further investigations of contact structures in the context of \emph{derived algebraic/symplectic geometry}.
   
    Regarding the study of shifted contact structures, we note that Maglio,  Tortorella and Vitagliano \cite{Maglio2024} have recently introduced and studied \emph{$0$-shifted} and \emph{$+1$-shifted contact structures} on \emph{differentiable stacks}, thus providing the foundations of {shifted contact geometry} in the stacky context.
     
     For shifted symplectic derived schemes, in particular, it has been  shown in \cite[Theorem 5.18]{Brav} that every $k$-shifted symplectic derived $\K$-scheme $(\bf X, \omega')$, with $k<0$, is Zariski locally equivalent to $(\spec A, \omega)$ for a pair $A,\omega$  in certain symplectic Darboux form \cite[Examples 5.8, 5.9 \& 5.10]{Brav}.

     In  \cite{BenBassat}, Ben-Bassat, Brav, Bussi, and Joyce extend the results of \cite{Brav} from  derived schemes to the case of derived Artin $\K$-stacks. 
     In this regard,  \cite[Theorem 2.8]{BenBassat} first proves that derived Artin $\K$-stacks also have nice atlases in terms of standard form cdgas. Then it is shown in \cite[Theorem 2.10]{BenBassat} that every $k$-shifted symplectic derived Artin $\K$-stack, with $k<0$, admits the so-called ``Darboux form atlas", and hence one has a Darboux-type theorem in this case as well. These results, with some side outcomes, motivate the current work.  
\paragraph{Main results and the outline.} In this work, with the same spirit as above, our goals are to extend the results of \cite{kib}  from derived schemes to derived Artin stacks and to provide several constructions of contact derived stacks. 

Recall that, in \cite{kib}, we introduced shifted contact structures on derived stacks and proved the following results, but \emph{only for}  (locally finitely presented) negatively shifted contact derived $\K$-schemes. In brief, we showed:

    \begin{theorem} \label{THM1}
    	Let $\bf X$ be a (locally finitely presented) derived $\K$-scheme.
    	\begin{enumerate}
    		\item[a.] \cite[Theorem 3.13]{kib}  Every $k$-shifted contact structure on $\bf X$, with $k<0$, is locally equivalent to $(\spec A, \alpha_0)$ for $A$ a minimal standard form cdga and $\alpha_0$ in a contact Darboux form.
    		\item[b.] \cite[Theorem 4.7]{kib}  Write $\mathcal{S}_{\bf{X}}$ for the total space  of a certain $\mathbb{G}_m$-bundle over  $\bf X$, constructed from the data of $k$-shifted contact structure on $\bf X$. Then  $\mathcal{S}_{\bf{X}}$ is a derived stack equipped with a $k$-shifted symplectic form $  \omega_{\bf{X}},$ which is canonically determined by the shifted contact structure of $\bf X$. We then call the pair $(\mathcal{S}_{\bf{X}}, \omega_{\bf{X}})$  the \emph{symplectification of $ \bf X $}.
    	\end{enumerate}
    \end{theorem}

In this paper, we  extend  Theorem \ref{THM1}   to the case of derived Artin $\K$-stacks 
as conjectured in \cite{kib}. In addition to that, we also provide several examples of contact derived stacks.
In short, the following theorems outline the main results of this paper:

\begin{theorem} \label{COROLLARY}
Theorem \ref{THM1}  also holds true for  negatively shifted contact derived Artin $\K$-stacks locally of finite presentation (cf. Theorems \ref{cor_darboux for Artin} \& \ref{cor_symplectization for Artin}).  
\end{theorem}

\begin{theorem}\label{ON EXAMPLES} Let $ \bf X $ be a derived Artin $\K$-stack locally of finite presentation. Denote by $ \mathbb{G}_a, \mathbb{G}_m $ the affine additive and multiplicative group schemes, respectively.
	\begin{enumerate}
		\item Let $\mathrm{T^*[n]} \bf X$ be the $n$-shifted cotangent stack. Then   the space $ \mathrm{J^1[n]{\bf X}:=T^*[n]}{\bf X }\times \mathbb{G}_a[n]$, called the \emph{$n$-shifted 1-jet stack} of $\bf X$, carries an $n$-shifted contact structure (cf. $\S$\ref{example_shifted jet space}).
		\item Let  $\pi_{{\bf X}}: \mathrm{T^*}{\bf X}\rightarrow \bf X$ be the natural projection. Given a prequantum 0-shifted Lagrangian fibration structure on $ \pi_{{\bf X}} $, there is a  $\mathbb{G}_m$-bundle on $ \mathrm{T^*}{\bf X}$ carrying a 0-shifted contact structure (cf. $\S$\ref{example_contact induced by prequantum on the cotangent}).
		\item Let  $\pi_{c_1(\mathcal{G})}: \mathrm{T^*}_{c_1(\mathcal{G})}\bf   X \rightarrow X$ be the $ c_1(\mathcal{G})$-twisted cotangent stack of $\bf X$, where $ c_1(\mathcal{G}) \in \mathcal{A}^1({\bf X}, 1) $ denotes the \emph{characteristic class} of a 0-gerbe $\mathcal{G}$ - a line bundle - on $\bf X$. Given a prequantum 0-shifted Lagrangian fibration structure on $ \pi_{c_1(\mathcal{G})} $, there is a  $\mathbb{G}_m$-bundle on $\mathrm{T^*}_{c_1(\mathcal{G})}\bf  X$ that carries a 0-shifted contact structure (cf. $\S$\ref{sec_example prequantum of twisted cotangent stacks}).
		\item Assume that $ G $ is a simple algebraic group
		over  $ \mathbb{K} $, and $ C $ be a  smooth and proper curve$/\mathbb{K} $. Consider the derived moduli stacks $LocSys_G(C), Bun_G(C)$ of \emph{flat $G$-connections on $C$}, \emph{principal $G$-bundles on $C$}. Then there is a  $\mathbb{G}_m$-bundle on $ LocSys_G(C) $ with a 0-shifted contact structure
		 (cf. $\S$\ref{sec_example stack of flat conn}).
	\end{enumerate}
\end{theorem}
Now, let us describe the content of this paper in more detail and provide an outline. In $\S$\ref{section_review of shifted symplectic strcs}, we review derived symplectic/contact geometries and the Darboux models (for derived schemes). 
In $\S$\ref{section_Artin stacs}, we will concentrate on  the stacky generalizations. In short, Section \ref{section_PTVV on Artin} outlines  the Darboux models for negatively shifted symplectic derived Artin stacks. With the same spirit, in Section \ref{section_Contact on Artin}, we consider the contact case and give the proof of Theorem \ref{COROLLARY} (cf. Theorems \ref{cor_darboux for Artin} \& \ref{cor_symplectization for Artin}). 
Finally, Section \ref{section_examples} provides several examples of contact derived stacks (with some background material on prequantization) and proves Theorem \ref{ON EXAMPLES}. We also have two appendices presenting some relevant constructions omitted in the main text.
    
    \paragraph {Acknowledgments.} I thank Alberto Cattaneo and  Ödül Tetik for helpful conversations. I am grateful to the Institute of Mathematics, University of Zurich, where an early draft of this paper was prepared during my stay there. %I personally benefited a lot from hospitality and research environment of the Institute. 
    The author acknowledges support of the Scientific and Technological Research Council of Turkey (T\"{U}BİTAK) under 2219-International Postdoctoral Research Fellowship (2021-1).
    
    Some parts of this paper were revisited and improved at the time when I joined Bilkent University Department of Mathematics in Fall 2023. Since then, I have benefited a lot from the discussions with Fırat Arıkan and Özgür Kişisel (both from Dept. of Math, Middle East Technical University). I am grateful to them for our research meetings.
    
    I wish to warmly thank the anonymous referee(s) for their valuable comments and suggestions, which helped a lot and improved the quality of the original manuscript.

    \paragraph{Conventions.} Throughout the paper, $ \mathbb{K} $ will be an algebraically closed field of characteristic zero. All cdgas will be graded in nonpositive degrees and  over $ \mathbb{K}.$ We assume that all classical $ \K $-schemes are \emph{locally of finite type}, and that all derived $ \K $-schemes/stacks $ \textbf{X} $ are  \emph{locally finitely presented.} %by which we mean that there exists a cover of $\bf X$ by affine opens $\spec  A$ with $A$ a finitely generated graded algebra. %I.e. $A$ has  finitely many generators and a finite number of relations over  $ \mathbb{K}.$
     
\section{Results for shifted  geometric structures on derived schemes} \label{section_review of shifted symplectic strcs} 
    
 \subsection{Some derived algebraic geometry} \label{section_prelim}
    We  provide a quick review on derived algebraic geometry. For details, we refer to \cite{Toen,ToenHAG,Luriethesis}. 
    
In this paper, we essentially use the functorial approach to define (higher) spaces of interest. It is very well known that using Yoneda's embedding,  spaces can be thought of as \textit{certain functors} in addition to the standard ringed-space formulation. In brief, we have the following diagram  \cite{Vezz2}: 
\begin{equation} \label{the diagram}
	\begin{tikzpicture}
	\matrix (m) [matrix of math nodes,row sep=1.5em,column sep=6.5em,minimum width=1.5 em] {
		CAlg_{\mathbb{K}}   & Sets  \\
		&  Grpds \\
		cdga_{\mathbb{K}} &  Ssets. \\};
	\path[-stealth]
	(m-1-1) edge  node [left] { } (m-3-1)
	edge  node [above] {{\small schemes}} (m-1-2)
	(m-1-1) edge  node [below] {} node [below] {{\small stacks}} (m-2-2)
	(m-1-1) edge  node [below] {} node [below] {{\small  higher stacks}} (m-3-2)
	(m-3-1) edge  node  [below] {{\small derived stacks}} (m-3-2)
	
	(m-1-2) edge  node [right] { } (m-2-2)
	(m-2-2) edge  node [right] { } (m-3-2);
	%edge [dashed,-] (m-2-1);
	\end{tikzpicture}
\end{equation} Here $ CAlg_{\mathbb{K}} $ denotes the category of commutative $\K$-algebras, and $ cdga_{\K}$ is the $\infty$-category\footnote{We actually mean the \emph{$\infty$-category} associated to the model category $ cdga_{\K} $,  with its natural model structure
for which equivalences are quasi-isomorphisms, and fibrations are
epimorphisms in strictly negative degrees.} of commutative differential graded $\mathbb{K}$-algebras in non-positive degrees. Denote by $St_{\K}$ the $\infty$-category of (higher) $\K$-stacks, where objects in $St_{\K}$ are defined via  Diagram \ref{the diagram} above. 

Recall that, in the underived setup, we  have the classical ``spectrum functor" 
\begin{equation*}
\cspec : (CAlg_{\mathbb{K}})^{op} \rightarrow St_{\K}.
\end{equation*} We then call an object $X$ of $ St_{\K} $ an \emph{affine $\K$-scheme} if $X\simeq \cspec A$ for some $A\in CAlg_{\mathbb{K}}$; and a \emph{$\K$-scheme} if it has an open cover by affine $\K$-schemes. In DAG, there also exists an appropriate concept of a \emph{spectrum functor}\footnote{In brief, it is the right adjoint to the global algebra of functions functor $\Gamma: dSt_{\K} \leftrightarrows cdga_{\K}^{op}: \spec$.} 
\begin{equation*}
\spec: cdga_{\K}^{op} \rightarrow dSt_{\K},
\end{equation*} which leads to the following definitions.
\newpage
\begin{definition}
	Denote by $dSt_{\mathbb{K}}$ the \bfem{$\infty$-category of derived stacks}, where  an object $ \textbf{X} $ of $dSt_{\mathbb{K}}$ is given as a certain $\infty$-functor\footnote{Using Yoneda's lemma, for a derived stack $\bf X$, we have  ${\bf X}: A\mapsto{\bf X}(A) \simeq Map_{dStk_{\mathbb{K}}}(\spec A, {\bf X})$, and hence any $A$-point $p\in {\bf X}(A) $ can be seen as a morphism $p: \spec A \rightarrow \bf X$ of derived stacks.} $ \textbf{X}:  cdga_{\mathbb{K}} \rightarrow  Ssets $ as in   Diagram \ref{the diagram}.
	More precisely, objects in $dSt_{\mathbb{K}}$ are \emph{simplicial presheaves on the site $(dAff)^{op}\simeq cdga_{\mathbb{K}}$  satisfying a descent condition}. For more details, we refer to \cite{ToenHAG}. 
\end{definition}

    \begin{definition}
    	An object $\textbf{X}$ in $dSt_{\mathbb{K}}$ is called an \bfem{affine derived $\mathbb{K}$-scheme} if  $\textbf{X}\simeq \spec A $ for some cdga $A \in cdga_{\K}$. An object $\bf X$ in $dSt_{\mathbb{K}}$ is then called a \bfem{derived $\mathbb{K}$-scheme} if it can be covered by Zariski open affine derived $\mathbb{K}$-schemes $Y \subset X$. Denote by $dSch_{\K} \subset dSt_{\mathbb{K}}$ the full $\infty$-subcategory of  derived $\mathbb{K}$-schemes, and we simply write $dAff_{\K} \subset dSch_{\mathbb{K}}$ for the full $\infty$-subcategory of  affine derived $\mathbb{K}$-schemes.% Note that $ \spec : (cdga_{\K}^{\infty})^{op} \rightarrow dAff_{\K} $ gives an equivalence of $\infty$-categories.
    \end{definition}

    %We should again note that throughout this paper, $ \mathbb{K} $ will be an algebraically closed field of characteristic zero. We also assume that all classical $ \K $-schemes are \emph{locally of finite type}, and all derived $ \K $-schemes $ \textbf{X} $ are  \emph{locally finitely presented}, by which we mean that $ \textbf{X} $ can be covered by Zariski open affines $\spec A,$ where $A$ is a finetely presented cdga over $\K.$

    %In addition to the spectrum functors $\spec, \cspec$ above, there is a natural \emph{truncation functor} $\tau: dSt_{\K} \rightarrow St_{\K}$, along with a fully faithfull left adjoint \emph{inclusion functor} $\iota: St_{\K} \hookrightarrow dSt_{\K}$, which can be thought of as an embedding of classical algebraic $\K$-spaces into derived spaces. In this regard, for a cdga $A$ there exists an equivalence $\tau \circ \spec A \simeq \cspec H^0(A)$. This means that if $\bf X$ is a (affine) derived $\K$-scheme, then its truncation $X=\tau(\textbf{X})$ is a (affine)  $\K$-scheme. Therefore, we can consider a derived $\K$-scheme $\bf X$ as \emph{an infinitesimal thickening of its truncation} $X$. It  follows that points of a derived $\K$-scheme $\bf X$ are the same as points of  of its truncation $X$. So the main difference between $X$ and $\bf X$ is in fact encoded by the scheme structure, not by the points!

\paragraph{Nice local models for derived $\K$-schemes.}  Let us first recall some basic concepts.
    
\begin{definition} \label{defn_standard forms}
	$ A \in cdga_{\mathbb{K}}$ is of \bfem{standard form} if  $A^0$ is a smooth finitely generated $\mathbb{K}$-algebra; the  module $\Omega^1_{A^0}$ of K\"{a}hler differentials is free $A^0$-module of finite rank; and the graded algebra $A$ is freely generated over $A^0$ by finitely many generators, all in negative degrees.

\end{definition}    
  
  In fact, there is a systematic way of constructing such cdgas starting from a smooth $\K$-algebra $A^0:=A(0)$ via the use of a sequence of localizations. More precisely, for any given $n\in \N$, we can inductively construct a sequence of cdgas 
  \begin{equation} \label{A(n) construction}
  A(0) \rightarrow A(1) \rightarrow \cdots \rightarrow A(i)\rightarrow \cdots A(n)=:A,
  \end{equation} where  $ A^0:=A(0) $, and $A(i)$ is obtained from $A(i-1)$ by adjoining generators in degree $-i$, given by $M^{-i}$, for all $i$. Here, each $M^{-i}$ is a free finite rank module (of degree $-i$ generators) over $A(i-1)$. Therefore, the underlying commutative graded algebra of $A=A(n)$ is freely generated over $A(0)$ by finitely many generators, all in negative degrees $-1,-2,\dots, -n$. For more details, we refer to \cite[Example 2.8]{Brav}.

 \begin{definition} \label{1st defn of minimalty}
 	A standard form cdga $ A $ is said to be \bfem{minimal} at  $p \in \cspec H^0(A)$ if $ A=A(n) $ is defined by using the minimal possible  numbers of graded generators in each degree $\leq 0$ compared to all other cdgas locally equivalent to  $ A $ near $p$. %(There will be an equivalent definition below, see Definition \ref{2nd defn of minimalty}.)
 \end{definition}  

 \begin{definition}
 	Let $A$ be a standard form cdga.  $A'\in cdga_{\mathbb{K}}$ is called a \bfem{localization} of $A$ if $A'$ is obtained from $A$ by inverting an element $f\in A^0$, by which we mean $A'=A\otimes_{A^0}A^0[f^{-1}]$. %for $f\in A^0$. 
 	
 	$A'$ is then of standard form with $A'^0 \simeq A^0[f]$. If $p\in \cspec H^0(A)$ with $ f(p)\neq0 $, we say $A'$ is a \emph{localization of $A$ at $p$}.
 \end{definition}
 With these definitions in hand, one has the following observations: 
 \begin{observation}
		Let $A$ be a standard form cdga. If $A'$ is  a \emph{localization} of $A$, then $\spec A' \subset \spec A$ is a Zariski open subset. Likewise, if $A'$ is a \emph{localization of $A$ at $p\in \cspec H^0(A)\simeq \tau(\spec A)$}, then  $\spec A' \subset \spec A$ is a Zariski open neighborhood of $p.$
	\end{observation}
\begin{observation}
	 Let  $ A$ be \emph{a standard form} cdga, then there exist generators $x_1^{-i}, x_2^{-i}, \cdots, x_{m_i}^{-i} $ in $A^{-i}$ (after localization, if necessary) with $i= 1, 2, \cdots, k$  \ and $m_i \in \mathbb{Z}_{\geq 0}$ such that 
	\begin{equation}
	A = A(0) \big[x_j^{-i} : i= 1, 2, \dots, k, \ j= 1,2, \dots, m_i\big],
	\end{equation} where the subscript $j$ in $x_j^i$ labels the generators, and the superscript $i$ indicates the degree of the corresponding element.  So, we can consider  $A$ as a \emph{graded polynomial algebra over $A(0)$ on finitely many generators, all in negative degrees.} 
\end{observation}%Moreover, since $A^0$ is finitely generated $\mathbb{K}$-algebra, there exists a surjection $\mathbb{K}[x_1^0, \cdots, x_{m_0}^0] \twoheadrightarrow A^0$ such that $A^0= \mathbb{K}[x_1^0,x_2^0, \cdots, x_{m_0}^0]/ I$ for some ideal $I$. %We shall be interested in the cases where the ideals $I$ are also finitely generated. Then $A^0$ is indeed said to be \emph{finitely presented}.
%	\item In general, we are interested in minimal standard form cdgas for which each number $m_i$ above is taken to be the least possible number such that $m_i=dim (H^i (\mathbb{L}_{A} |_p))$  and $d^i |_p =0$ for $i= -1, -2, \cdots, k$, and $p \in spec (H^0(A))$. 

%\begin{definition}
%	We then define the \emph{virtual dimension} of $A$ to be the integer $\vdim A= \sum _i (-1)^im_i$. 
%\end{definition}
	
The following theorem outlines the central results from \cite[Theorem 4.1 \& 4.2]{Brav} concerning the construction of useful local algebraic models for  derived $\K$-schemes.  The upshot is that given a derived $\K$-scheme $\bf X$ (locally of finite presentation) and a point $x\in \bf X$, one can always find a ``refined"  affine neighborhood $ \spec A $ of $x$, which is very useful for explicit presentations.    
In short, we have:
\begin{theorem} \label{localmodelthm}
	Every derived $\mathbb{K}$-scheme $X$ is Zariski locally modelled on $ \spec A $ for a minimal standard form cdga $ A$. 
\end{theorem}

\paragraph{Nice local models for cotangent complexes of derived schemes.} Given	$ A \in cdga_{\mathbb{K}}$, $ d$ on $A$ induces a differential on $\Omega_A^1$, denoted again by $ d$. This makes $\Omega_A^1$ into a $\rm dg$-module $(\Omega_A^1, d)$ with the property that $\delta \circ  d = d \circ \delta,$ where $\delta: A \rightarrow \Omega_A^1$ is the universal derivation of degree zero. Write  the decomposition of $ \Omega^1_{A} $ into graded pieces 
$ \Omega^1_{A} =  \bigoplus_{k=-\infty}^0 \big(\Omega^1_{A}\big)^k $
with the  differential $d: \big(\Omega^1_{A}\big)^k \longrightarrow \big(\Omega^1_{A}\big)^{k+1}$. Then we define the \textit{de Rham algebra of $A$} as a  double complex 
   \begin{equation}
   DR(A)= Sym_A(\Omega_A^1[1]) \simeq \bigoplus \limits_{p=0}^{\infty} \bigoplus \limits_{k=-\infty}^0 \big(\Lambda^p \Omega^1_{A}\big)^k [p],
   \end{equation} where the gradings $p, k$ are called the \emph{weight} and the \emph{degree}, respectively. Also, there are two derivations (differentials)  on $DR(A)$, namely the \emph{internal differential d}$ : \big(\Lambda^p \Omega^1_{A}\big)^k [p] \longrightarrow \big(\Lambda^p \Omega^1_{A}\big)^{k+1} [p] $ and the \emph{de Rham differential} $d_{dR}:\big(\Lambda^p \Omega^1_{A}\big)^k [p] \longrightarrow \big(\Lambda^{p+1} \Omega^1_{A}\big)^k [p+1]$ such that $d_{tot}= d + d_{dR}$ and  \begin{equation}
   d^2=d_{dR}^2=0, \text{\ and \ } d\circ d_{dR} =- d_{dR} \circ d.
   \end{equation} Here, one also has the natural multiplication on $ DR(A)$:
   \begin{equation}
   \big(\Lambda^p \Omega^1_{A}\big)^k [p] \times \big(\Lambda^q \Omega^1_{A}\big)^{\ell} [q] \longrightarrow \big(\Lambda^{p+q} \Omega^1_{A}\big)^{k+\ell} [p+q].
   \end{equation} 
   
 % \begin{observation}
 % 	 It should be noted that the constructions of $\Omega_A^1, \ DR(A)$ depend only on the underlying commutative graded algebra of $A$, not on the differential $d$ on $A.$ 
  %\end{observation}

   Note that even if both $\mathbb{L}_{A} \text{ and } \Omega^1_{A}$ are closely related, the identification of $ \mathbb{L}_{A} $ with  $ \Omega^1_{A} $ is not true for an arbitrary   $ A \in cdga_{\mathbb{K}}$ \cite{Brav}. But, when $A=A(n)$ is a standard form cdga, we have the following description for the restriction of the cotangent complex $\mathbb{L}_A $ to $\cspec H^0(A)$. In this paper, we only give a brief version. More details and the proof can be found in \cite[Prop. 2.12]{Brav}.

 \begin{proposition}  \label{proposition_L as a complex of H^0 modules}
If $A=A(n)$ with $n\in \N$ is a standard form cdga constructed inductively as in (\ref{A(n) construction}), then the restriction of $\mathbb{L}_A $ to $\cspec H^0(A)$ is represented by a  complex of free $H^0(A)$-modules.

 \end{proposition}  

%\begin{definition} \label{2nd defn of minimalty}
%Let $A=A(n)$ with $n\in \N$ be a standard form cdga constructed inductively. $A$ is said to be \emph{minimal at} $p\in \cspec H^0(A)$ if the internal differential $d^{-i}|_p = 0$ in the complex representing $\mathbb{L}|_{\cspec H^0(A)}$.
%\end{definition}

    %Note that Definition \ref{2nd defn of minimalty} implies  $m_i=\dim (H^i (\mathbb{L}_{A} |_p))$ for each $i$, and hence $ A $ is defined by using the minimum number of graded variables in each degree $\leq 0$ compared to all other cdgas locally equivalent to  $ A $ near $p$. Therefore, one can recover Definition \ref{1st defn of minimalty}.

\subsection{Derived symplectic geometry and local models}   \label{secion_PTVV's symplectic geometry on spec A} 
  %Let $\bf X$ be a locally finitely presented derived $\K$-scheme with $p\geq 0, \ k\in \Z$. 
  Pantev et al. \cite{PTVV} define the \emph{simplicial sets of $p$-forms of degree $k$  and closed $p$-forms of degree $k$}  on derived stacks. Denote these simplicial sets by $\mathcal{A}^p({\bf X},k)$ and $\mathcal{A}^{p,cl}({\bf X},k)$, respectively. These definitions are in fact given first for affine derived $\K$-schemes. Later, both concepts are defined for  general derived stacks $\bf X$ in terms of mapping stacks $\mathcal{A}^p(\bullet,k)$ and $\mathcal{A}^{p,cl}(\bullet,k)$, respectively.% A summary of key ideas can be found in \cite[Section 3.4]{Brav}
  
 The space $\mathcal{A}^{p,cl}({\bf X},n)$ of closed $ p $-forms on a general derived stack $\bf X$ can be a rather complicated even when $\bf X$ is a nice derived Artin stack. However, \cite[Prop.  1.14]{PTVV} gives the following identification for the space $ \mathcal{A}^p({\bf X},n) $ of $n$-shifted $ p $-forms: \begin{equation*}
 	\mathcal{A}^p(X,n)\simeq Map_{QCoh(X)}(\mathcal{O}_X, \wedge^p\mathbb{L}_X[n]). 
 \end{equation*}
  
  Let $ \textbf{X} = \spec A$ with $A$ a standard form cdga\footnote{It should be noted that the results that are cited or to be proven in this section are all about the \emph{local structure} of derived schemes. Thus, it is enough to consider the (refined) affine case.}, then take $ \Lambda^p\mathbb{L}_A=\Lambda^p \Omega^1_{A}$.  Therefore,  elements of $\mathcal{A}^p({\bf X},k)$  form a simplicial set such that $k$-cohomology classes of the complex $ \big(\Lambda^p \Omega^1_{A}, d\big) $  correspond to the connected components of this simplicial set. Likewise,  the connected components of $\mathcal{A}^{p,cl}({\bf X},k)$ are identified with the $k$-cohomology classes of the complex $\prod_{i\geq0} \big(\Lambda^{p+i} \Omega^1_{A}[-i], d_{tot}\big)$. Then we have: %the following definitions. %We want to work with explicit representatives for these cohomology classes. 

 \begin{definition} \label{defn_p form of deg k}
 	Let $\textbf{X}=\spec A$ be an affine derived $\mathbb{K}$-scheme for $A$ a minimal standard form cdga. A \bfem{$k$-shifted $p$-form on $\bf X$} for $p\geq 0$ and $k \leq 0$ is an element 
 	$ \omega^0 \in \big(\Lambda^p \Omega^1_{A}\big)^k \text{ with } d\omega^0=0.  $
 	
 \end{definition} Note that an element  $\omega^0$ defines a cohomology class $[\omega^0] \in H^k\big(\Lambda^p \Omega^1_{A}, d\big),$ where two $ p $-forms  $\omega_1^0, \omega_2^0 $ of degrees $k$ are \emph{equivalent} if $\exists \alpha_{1,2} \in \big(\Lambda^p \Omega^1_{A}\big)^{k-1}$ s.t. $\omega_1^0-\omega_2^0= d\alpha_{1,2}$.

   \begin{definition}\label{defn_closed p form}
   Let $\textbf{X}=\spec A$ be an affine derived $\mathbb{K}$-scheme with $A$ a minimal standard form cdga. A \bfem{closed $k$-shifted $p$-form on $\bf X$} for $p\geq 0$ and $k \leq 0$ is a sequence $\omega=(\omega^0, \omega^1, \cdots)$ with $\omega^i \in \big(\Lambda^{p+i} \Omega^1_{A}\big)^{k-i}$ such that $d_{tot}\omega=0$, which splits according to weights as $d\omega^0=0$ in $ \big(\Lambda^p \Omega^1_{A}\big)^{k+1}$
   and $d_{dR}\omega^i + d\omega^{i+1}=0$ in $\big(\Lambda^{p+i+1} \Omega^1_{A}\big)^{k-i}$, $i \geq 0.$ 
   
   That is, a closed $k$-shifted $p$-form consists of an actual $k$-shifted $p$-form $ \omega^0 $ and the data $ (\omega^i)_{i>0} $ of $ \omega^0 $ being coherently $d_{dR}$-closed.
  It then follows that there also exists a natural projection morphism  
 	$ \pi: \mathcal{A}^{p,cl}(\textbf{X},k)\longrightarrow \mathcal{A}^{p}(\textbf{X},k), \ \ \omega=(\omega^i)_{i\geq 0} \longmapsto \omega^0. $

   \end{definition}

\begin{definition}
	A closed $k$-shifted $2$-form $ \omega=(\omega^i)_{i\geq 0} $ on  $\textbf{X}=\spec A$ for a (minimal) standard form cdga $A$ is called a \bfem{$k$-shifted symplectic structure} if  the induced map $$ \omega^0\cdot: \mathbb{T}_{A} \rightarrow \Omega^1_{A}[k], \ Y \mapsto \iota_{Y} \omega^0,$$ is a quasi-isomorphism, where $\mathbb{T}_{A}=(\mathbb{L}_{A})^{\vee}\simeq\Hom_{A}(\Omega^1_{A},A)$\footnote{Thanks to the identification $\mathbb{L}_A\simeq (\Omega_A^1, d)$ for $A$ a (minimal) standard form cdga.} is the \textit{tangent complex of $A$.} The requirement for the induced map $ \omega^0\cdot $ is called the \bfem{non-degeneracy condition}.
\end{definition}

\paragraph{Symplectic Darboux models for derived schemes.} \label{the pair}
One of the main theorems in  \cite{Brav} provides a $ k $-shifted version of the classical Darboux theorem in symplectic geometry. The statement is as follows. 

\begin{theorem} (\cite[Theorem 5.18]{Brav}) \label{Symplectic darboux}
	Given a derived $\mathbb{K}$-scheme $\bf{X}$ with a $k$-shifted symplectic form $\omega'$ for $k<0$ and $x\in \bf X$, there is a local model $\big(A, f: \spec A \hookrightarrow \bf{X}, \omega \big)$  and $p \in \cspec (H^0(A))$ such that $f$ is an open inclusion with $f(p)=x$, $A$ is a standard form that is minimal at  $p$, and $\omega$ is a $k$-shifted symplectic form on $\spec A$ such that $A, \omega$ are in Darboux form, and $f^*(\omega') \sim \omega$  in $ \mathcal{A}^{2,cl}(\textbf{X},k).$
\end{theorem} 

In fact, Theorem \ref{Symplectic darboux} shows that such $ \omega $  can be constructed explicitly depending on the integer $k<0$. Indeed, there are three cases in total: 
  $(1) \ k  $  is odd;  $ (2) \ k/2 $  is even;  and  $ (3) \ k/2  $ is odd. 
 For instance, when $k$ is odd, one can find a minimal standard form cdga $A$, with $ `` $coordinates" $x^{-i}_j, y^{k+i}_j \in A$, and a Zariski open inclusion $f: \spec A \hookrightarrow \bf X$ so that $f^*(\omega') \sim \omega=(\omega^0, 0, 0, \dots)$, where
$  \omega^0= \sum_{i,j} d_{dR}x_j^{-i} d_{dR}y_j^{k+i}.$ We will not give any further detail on the aforementioned cases in this paper. Instead, we refer to  \cite[Examples 5.8, 5.9, and 5.10]{Brav}.

We just wish to present a  result  that plays a significant role in constructing Darboux-type local models. The upshot is that one can always simplify the  given closed 2-form $ \omega= (\omega^0,\omega^1, \omega^2, \dots)$ of degree $k<0$ on $\spec A$ so that $\omega^0$ can be taken to be exact and $\omega^i=0$ for all $i>0.$ More precisely, we have\footnote{The result is based on the interpretation of such forms in the context of \emph{cyclic homology theory of mixed complexes} \cite{PTVV,Brav}.}:

\begin{proposition} (\cite[Prop. 5.7]{Brav} )\label{Proposition_exactness}
	Let  $ \omega= (\omega^0,\omega^1, \omega^2, \dots)$ be a closed 2-form of degree $k<0$ on $\spec A$ for $A$ a standard form cdga over $\K.$ Then there exist $H\in A^{k+1}$ and $\phi \in (\Omega^1_{A})^k$ such that $dH=0$ in $A^{k+2}$, $d_{dR}H+d\phi=0$ in $(\Omega^1_{A})^{k+1}$, and $\omega \sim (d_{dR}\phi, 0, 0, \dots).$ %In fact, $d_{dR}\phi=k\omega^0.$ 
	
	%Moreover, if $(H', \phi')$ is another such pair for fixed $\omega, k, A,$ then there exist $h\in A^k$ and $\sigma\in (\Omega_A^1)^{k-1}$ such that $H-H'=dh$ and $\phi - \phi'=d_{dR}h+d\sigma.$
\end{proposition}

\subsection{Derived contact geometry and local models} \label{section_shifted contact structures and Darboux forms}
In this section, we review the central constructions and results from \cite{kib}. In a nutshell, a \textit{ $k$-shifted contact structure} on a derived Artin stack  consists of a  morphism $f:\mathcal{K}\rightarrow \mathbb{T}_{{\bf X}}$ of perfect complexes, a line bundle $L$ such that $Cone(f)\simeq L[k]$, and a locally defined $k$-shifted 1-form $\alpha$ satisfying a non-degeneracy condition. With this structure in hand, \cite{kib} presents a Darboux-type theorem and the  process of canonical symplectification for \emph{negatively shifted contact derived schemes}. 
\subsubsection{Basic concepts}
 %More precisely, we have:

%\paragraph{Preliminaries.} 
Let us start with some terminology. Recall that the definition of  a \emph{derived Artin $\K$-stack} can be formalized by using an inductive construction \cite[Section 5.1]{Luriethesis}. Roughly speaking, we call ${\bf X} \in dSt_{\K}$ a \emph{derived Artin $\K$-stack} if it can be locally represented by an affine derived $\K$-scheme with respect to the ``smooth topology". Thus, we require the existence of a ``smooth" surjection $\varphi: U \rightarrow \bf X$ (of some relative dimension $m$) so that $U$ is a disjoint union of affine derived schemes. 

In order to the make sense of the smoothness of $\varphi$ as above, we require that the fibers of the morphism $\varphi$ are already suitable geometric objects. To this end, one should start with some subcategory $S_0$, called \emph{$0$-stacks}, and then define $S_{n+1}$ to be the class of objects $\bf X$ in $ dSt_{\K}  $ having a ``smooth" surjection $\varphi: U \rightarrow \bf X$ with $U$  a disjoint union of affine derived schemes such that each fiber $U \times_{\bf X} \spec A$ lies in the class $S_n.$ 

 Technically speaking, 	an object $X \in dSt_{\K}$ is called a \emph{derived Artin $\K$-stack} if it is $m$-geometric for some $m$, and the underlying classical stack is 1-truncated (i.e. it is just a stack, not higher stack). For  details, we refer to \cite[$\S$ 5.1]{Luriethesis} or \cite[$\S$ 1.3.3]{ToenHAG}.% More formally, we have

The upshot is that any such object $\bf X$ of $ dSt_{\K}$ comes with a smooth surjective morphism  $\varphi: U \rightarrow \bf X$ with $U$ a  derived $\K$-scheme. We call such morphism an \bfem{atlas.} Therefore, the following definition will be sufficient for our purposes.

%\begin{definition} %(\cite[Def. 1.3.3.1]{ToenHAG} or \cite[Section 5.1]{Luriethesis})
%	An object $X \in dSt_{\K}$ is called a \textit{derived Artin $\K$-stack} if it is $m$-geometric for some $m$.
%\end{definition} 

\begin{definition}
	By a \bfem{derived Artin $\K$-stack},  we mean an object $\bf X$ of $ dSt_{\K}$ possessing an atlas (smooth of some relative dimension) near each point of $\bf X$.
\end{definition}

Now, let us introduce \emph{shifted contact structures} on derived Artin stacks:

\begin{definition}\label{defn_preshiftedcontact} Let $\bf X$ be a locally finitely presented derived (Artin) stack.
	A \bfem{pre-$ k $-shifted contact structure} on $\bf X$ is given by  a shifted line bundle $L[k]$ with a morphism $\alpha:\mathbb{T}_{{\bf X}}\rightarrow L[k]$. Denote such a structure by $(L[k],\alpha)$.%a perfect complex $\mathcal{K}$ on $\bf X$ with a monomorphism $\kappa: \mathcal{K} \rightarrow \mathbb{T}_{{\bf X}}$ of perfect complexes whose cone $Cone(\kappa)$ is weak equivalent to $L[k]$, where $L$ is a line bundle. Denote such a structure on $\bf X$ simply by $(\mathcal{K}, \kappa, L)$.
\end{definition}

%\begin{notation}
Note that we can consider a pre-$ k $-shifted contact data as a perfect complex $\mathcal{K}$ and a line bundle $L$ along with a morphism $\kappa: \mathcal{K}\rightarrow\mathbb{T}_{{\bf X}}$   such that $Cone(\kappa)\simeq L[k]$. Here, we have  a  cofiber sequence $\mathcal{K} \rightarrow \mathbb{T}_{{\bf X}}\rightarrow L[k]$ in $QCoh({\bf X})$. Since $QCoh({\bf X})$ is a stable $\infty$-category, the cocone of $ \mathbb{T}_{{\bf X}}\rightarrow L[k]$ is equivalent to $\mathcal{K}$. 
We then may denote  a pre-$ k $-shifted contact structure on $\bf X$  by $(\mathcal{K}, \kappa, L)$.
%\end{notation}
\begin{definition}\label{defn_shiftedcontact}
	We say that a pre-$ k $-shifted contact structure $(\mathcal{K}, \kappa, L)$ on $\bf X$ is a \bfem{$ k $-shifted contact structure} if locally on $\bf X$, where $L$ is trivial, the induced $k$-shifted 1-form\footnote{We can locally identify the map $\alpha$ with the induced shifted one-form using the  trivialization of $L^{\vee}[k]$.} $\alpha: \mathbb{T}_{{\bf X}} \rightarrow \mathcal{O}_{{\bf X}}[k]$ is such that the map $d_{dR}\alpha|_{\mathcal{K}} \ \cdot:= \kappa^{\vee}[k] \circ( d_{dR}\alpha \ \cdot) \circ \kappa : \mathcal{K}\rightarrow\mathcal{K}^{\vee}[k]$ is a weak equivalence. 
	
	In that case, we say the $k$-shifted 2-form $d_{dR}\alpha$ is \bfem{non-degenerate on $\mathcal{K}$}. Also, we call such local form a \bfem{$k$-contact form}.
\end{definition}

\begin{remark}
	When $k\leq 0$, the triangle $\mathcal{K} \rightarrow \mathbb{T}_{{\bf X}} \rightarrow L[k]$ splits locally for any affine derived scheme (so, this also holds Zariski
	locally for any derived scheme)\footnote{We thank the anonymous referee for this remark.}. In fact,  the nondegeneracy condition implies that $\mathcal{K}$ has $\mathrm{Tor}$-amplitude
	$ [0,-k] $ so that $\mathcal{K}[-k]$ is connective. Then the connecting homomorphism $ L[k] \rightarrow  \mathcal{K}[1] $ in the exact triangle is equivalently $ L \rightarrow  \mathcal{K}[1-k] $. Notice that $ \mathcal{K}[1-k] $ is concentrated in degrees $\leq-1$, so this morphism
	is automatically zero on any  affine derived scheme, which implies the desired splitting.
\end{remark}
Let $\bf X$ be a locally finitely presented  derived Artin stack with a $k$-shifted contact structure $(\mathcal{K}, \kappa, L)$. Recall from Yoneda's lemma, ${\bf X}(A) \simeq Map_{dPstk}(\spec A, {\bf X})$, and hence any $A$-point $p\in {\bf X}(A) $ can be seen as a morphism $p: \spec A \rightarrow \bf X$ of derived pre-stacks. Then, let us consider the pair $ (p, \alpha_p )$, with $  p \in {\bf X} (A), \ \alpha_p \in p^*(\mathbb{L}_{{\bf X}}[k]), $ such that  $ Cocone(\alpha_p)\simeq \mathcal{K}$. For $A \in cdga_{\K}$, there is a $\mathbb{G}_m(A)$-action on the pair $ (p, \alpha_p )$ by 
\begin{equation*}
	f \triangleleft (p, \alpha_p):= (p, f \cdot \alpha_p).
\end{equation*} Denote by $H^0$  the functor sending $A \mapsto H^0(A).$ Denote the image under $H^0$ of an element $f$ simply by $f^0.$ Note that localizing $ A $ if necessary, w.l.o.g. we may assume that the image $f^0$ is always invertible. It follows that $f^0$ lies in $(A^0)^{\times}$, which is by definition $\mathbb{G}_m(A^0)=(A^0)^{\times}$.

\begin{observation} \label{observation_cocones}
	If ${\bf X}, (p, \alpha_p )$, and the $\mathbb{G}_m(A)$-action are as above, then for an element $f\in \mathbb{G}_m(A)$, we can obtain $Cocone(f\cdot \alpha_p) \simeq Cocone(\alpha_p)$ by using the invertibility of $f.$ 
\end{observation}

From Proposition \ref{proposition_L as a complex of H^0 modules}, on a refined affine neighborhood, say $\spec A$ with $A$ a minimal standard form cdga, the perfect complexes $\mathbb{T}_A, \mathbb{L}_A$, when restricted to ${\cspec H^0(A)}$, are both free finite complexes of $H^0(A)$-modules. In that case, Definitions \ref{defn_preshiftedcontact} and \ref{defn_shiftedcontact}, and Observation \ref{observation_cocones} will reduce to the following local descriptions, where  $\mathcal{K}$ is now just equivalent to the usual $\ker \alpha$ in $Mod_A$; and  $L$ in the splitting corresponds to the  line bundle generated by the Reeb vector field of the classical case.  

More precisely, from \cite[$\S 3.2$]{kib}, when restricted to the (nice) local models, we equivalently have the following proposition/definition:
\begin{proposition} \label{defn_contact strc on good affines}
	\emph{(Shifted contact structures with nice affine models)} For a (minimal) standard form cdga $A$ and $k<0$, any $k$-shifted contact structure on  $\textbf{X}=\spec A$ can be strictified in the sense that the resulting contact data consists of \begin{itemize}
		\item a submodule $\mathcal{K}$ with the natural inclusion $i:\mathcal{K}\hookrightarrow Der(A)$
 such that $Cone(i)\simeq coker(i)$ is the quotient complex and of the form $L[k]$, with $L$ a line bundle; and
 \item  a  $k$-shifted 1-form $\alpha$  on $\spec A$  with the property that $\mathcal{K}\simeq \ker \alpha$ so that the $k$-shifted 2-form $d_{dR}\alpha$ is non-degenerate  on $\ker \alpha$.
	\end{itemize}Here $Der(A)=(\Omega^1_{A})^{\vee}=\Hom_{A}(\Omega^1_{A},A)$, where $ \Omega^1_{A} $ is the $A$-module K\"{a}hler differentials such that $ \Omega^1_{A}|_{\cspec H^0(A)} $ is represented as a (bounded) complex of free $H^0(A)$-modules by Proposition \ref{proposition_L as a complex of H^0 modules}. 
In that case, over $ p\in \cspec H^0(A) $, one has the natural splitting 	
$ 	Der(A)|_{\cspec H^0(A)}=  \ker \alpha |_{\cspec H^0(A)} \oplus L[k]|_{\cspec H^0(A)}. $ 	
Adopting the classical terminology, we sometimes call the sub-module $\ker \alpha$  above a \bfem{$ k $-shifted (strict) contact structure with the defining $ k $-contact form}  $ \alpha $.
\end{proposition}
%In what follows, we give a prototype construction for  $k$-shifted contact forms.
% as in the case of shifted symplectic Darboux-like local models.
\begin{example} \emph{($k$-contact forms, with odd $k<0$)}\label{model example}
	In this example, fixing $\ell\in \N$, we will construct a (minimal) standard form cdga carrying a $k$-shifted contact structure  for $k=-2\ell -1$. In brief, we will coherently extend the symplectic case \cite[Example 5.8]{Brav}. To this end, we make use of similar notations and constructions from that example.    Modifications for even shifts are outlined in Appendix \ref{app. non-odd shifts}. 
	
\bfem{Step-1: Construction of a ``contact" cdga.}	Let $ A(0) $ be a smooth $ \K $-algebra  of dimension $ m_0 $. Assume that there exist degree 0 variables $x_1^0, x_2^0, \dots, x_{m_0}^0$ in $A(0)$ defining global \'{e}tale coordinates $(x_1^0, x_2^0, \dots, x_{m_0}^0): \spec A(0) \rightarrow \mathbb{A}^{m_0}$ on $\spec A(0)$ such that $d_{dR}x_1^0,\dots,d_{dR}x_{m_0}^0$ form a $A(0)$-basis for $\Omega_{A(0)}^1$.  
	Next, choosing non-negative integers $m_1,\dots, m_{\ell}$, we define   a \bfem{commutative graded algebra} $A$ to be the free graded $\K$-algebra over $A(0)$ generated by the variables  
	\begin{align} 
	& x_1^{-i}, x_2^{-i}, \dots, x_{m_i}^{-i}& &\text{ in degree } -i \ \ \ \text{ for } i= 1,  \dots, \ell, \label{var set1} \\
	& y_1^{k+i}, y_2^{k+i}, \dots, y_{m_i}^{k+i}& & \text{ in degree } k+i \ \text{ for } i=1,\dots, \ell, \label{var set2} \\
	& z^k, y_1^{k}, y_2^{k}, \dots, y_{m_0}^{k}& & \text{ in degree } k,
	\end{align} such that  $\Omega^1_{A}$ is  the free  $A$-module of finite rank with basis $ \{d_{dR}x_j^{-i}, d_{dR}y_j^{k+i}, \dR z^k: \forall i,j\}.$ Here, we call $z^k$  the \bfem{distinguished variable} (of $\deg k)$. Also, we choose an element $H\in A^{k+1}$  satisfying  the \emph{classical master equation} \begin{equation} \label{defn_CME}
	\displaystyle \sum_{i=1}^{\ell} \sum_{j=1}^{m_i} \dfrac{\partial H}{\partial x_j^{-i}} \dfrac{\partial H}{\partial y_j^{k+i}}=0 \text{ in } A^{k+2}. 
\end{equation} Here, due to  degree reasons, $H$ does not involve  any of $z^k,y^k_j$'s.  We call such $H$ the \emph{Hamiltonian}.
Then we define the \bfem{internal differential} $d$ on $A$ by  the equations \begin{align} \label{defn_internal d contact}
	d|_{A(0)}&=0; \ dx_j^{-i} =  \dfrac{\partial H}{\partial y_j^{k+i}} \text{ for all } i>0,j; \  \ dy_j^{k+i} =  \dfrac{\partial H}{\partial x_j^{-i}} \text{ for all } i,j; \text{ and } \nonumber \\ -kdz^k&= H+d\Big[\sum_{i,j} (-1)^{i} ix_j^{-i} y_j^{k+i} \Big].
\end{align} Notice that the condition on $H$ implies  $d^2=0$ on each generator \cite{Brav}. Also, we have  $\vdim(A)=-1$.

\bfem{Step-2: Pre-contact data.} Next, we introduce the element $\alpha \in (\Omega_{A}^1)^k$ given by
\begin{equation} \label{the local primitive contact model}
	\alpha= d_{dR}z^k+ \displaystyle \sum_{i=0}^{\ell} \sum_{j=1}^{m_i} y_j^{k+i}d_{dR}x_j^{-i}.
\end{equation} Let us first show that $d\alpha=0$.   To this end, we compute $d_{dR}\alpha=\sum_{i=0}^{\ell} \sum_{j=1}^{m_i} d_{dR}x_j^{-i}d_{dR}y_j^{k+i}$. Then  the element   $d_{dR}\alpha$ is  
$(d+\dR)$-closed  by \cite[Example 5.8]{Brav}.  The same example also implies that  if we let\footnote{Alternatives like $ k\sum_{i,j} y_j^{k+i}d_{dR}x_j^{-i} $ can be obtained by replacing $H, \phi$ by suitable $H+d[\cdots]$ and $\phi+\dR [\cdots]$.} $$\phi:=  \sum_{i=0}^{\ell} \sum_{j=1}^{m_i} \big[(-i)x_j^{-i}d_{dR}y_j^{k+i}+(k+i)y_j^{k+i}d_{dR}x_j^{-i}\big]$$ and $ H $ as above, then the pair   $(\phi,H)\in (\Omega^1_A)^k\times A^{k+1}$ is a solution to the equations \begin{equation} \label{important relations}
dH=0 \text{ in } A^{k+2}, \ d_{dR}H+d\phi=0 \text{ in } (\Omega^1_{A})^{k+1}, \text{ and } d_{dR}\phi=kd_{dR}\alpha.
\end{equation} Observe that   $  \phi  + \dR \big[ \sum_{i,j} (-1)^{i} ix_j^{-i} y_j^{k+i} \big] =k\sum_{i,j} y_j^{k+i}d_{dR}x_j^{-i}$, then we can write 
\begin{equation} \label{local contact model}
	k\alpha=kd_{dR}z^k+ \displaystyle k\sum_{i=0}^{\ell} \sum_{j=1}^{m_i} y_j^{k+i}d_{dR}x_j^{-i}=k\dR z^k + \phi  + \dR \Big[ \sum_{i,j} (-1)^{i} ix_j^{-i} y_j^{k+i} \Big].
\end{equation}  %Here $k\alpha \in \Omega^1_{A}[k]$ is a representative of $\phi$, but with  $d\phi=-d_{dR} H$. 
 As we set ${-kdz^k=H+d[\cdots]}$ in (\ref{defn_internal d contact}), we obtain $d(k\alpha)=\dR \circ d(-kz^k) + d\phi + d \circ \dR [\cdots] = \dR(H+d[\cdots])-\dR H - \dR \circ d [\cdots]=0$ using (\ref{important relations}). In that case, $\alpha$ is then \emph{$d$-closed, and hence a 1-form of degree $k.$} In other words, the pair $(H,k\alpha-kd_{dR}z^k-\dR[\cdots])$, with $H,\alpha$ as above, provides a solution to the  equations (\ref{important relations}), which implies  $d\alpha=0$.
%That is,  we have a pair of solution $(H,\phi)$ to the defining equations (\ref{important relations}) such that $H$ is an element in $A^{k+1}$ with   and ${\phi=\alpha-kd_{dR} z^k}$.

Next, we show that there is a canonical pre-$k$-shifted contact structure on $\spec A$ induced by $\ker \alpha.$ For $i=0,\dots,\ell, \text{ and } 1\leq j\leq m_i$,
denote the vector fields annihilating $\alpha$ by
\begin{align}
	%\sigma_{ij}^{0} &=ix^{-i}_j\partial/\partial x_j^{-i} + (k+i)y^{k+i}_j\partial/\partial y^{k+i}_j & & \text{ in degree }  0, \nonumber \\
	\zeta_j^{i}&= \partial/\partial x_j^{-i}-y^{k+i}_j\partial/\partial z^{k} & & \text{ in degree } i, \nonumber \\
	\eta^{-k-i}_j &= \partial/\partial y^{k+i}_j & & \text{ in degree }  -k-i.
\end{align} %Observe that $ \sigma_{ij}^{0}$ can  be written by using the other vector fields\footnote{Observe that $ k\sigma_{ij}^{0}= ix^{-i}_j	\zeta_j^{i}+ (k+i)y^{k+i}_j\eta^{-k-i}_j$.} $ 	\zeta_j^{i} \text{ and }  \eta^{-k-i}_j$. 
We thus obtain $\ker \alpha= \Span_{A} \{ \zeta_{j}^{i}, \ \eta^{-k-i}_{j} \} \hookrightarrow  \mathbb{T}_A $ over $\spec A$, with the  $\mathrm{Tor}$-amplitude $[0,-k]$. From Prop. \ref{proposition_L as a complex of H^0 modules}, the restrictions of ${\mathbb{T}_A},{\ker \alpha}$ to ${\spec H^0(A)}$ are both complexes of free $H^0(A)$-modules. Likewise, the restriction of $ \ \faktor{\mathbb{T}_A}{\ker \alpha} \ $ to ${\spec H^0(A)}$ is then generated by the $\deg -k$ vector field $\partial/\partial z^{k}$ as an $H^0(A)$-module, and hence it is equivalent to the complex concentrated in degree $-k$. It follows that the quotient  can also be identified with a $k$-shifted line bundle $ L[k] $ on  $\spec A$ such that  we get a  fiber-cofiber sequence \begin{equation*}
\ker \alpha \rightarrow \mathbb{T}_A \rightarrow L[k],
\end{equation*}and hence a map  $ \mathbb{T}_A \rightarrow L[k].$ 
Therefore, the data of $L[k]$ and the map $ \mathbb{T}_A \rightarrow L[k]$ together with the $k$-shifted 1-form $\alpha$ above define a \bfem{pre-$k$-shifted contact structure} on $\spec A$  in the sense of Definition \ref{defn_preshiftedcontact}.  
Call such $\alpha$  a \bfem{pre-$k$-contact form}.

\bfem{Step-3: From pre-contact to contact.} So far, we have obtained a pre-$k$-contact form $\alpha \in (\Omega_A^1)^k$ as in (\ref{the local primitive contact model}) using  the  variables in (\ref{var set1}), the differential $d$ defined by (\ref{defn_internal d contact}), and the  equations in (\ref{important relations}). It remains to show that $d_{dR}\alpha$ { is non-degenerate on }  $ {\ker \alpha} $.
To this end,  it suffices to prove the non-degeneracy of the induced map  $d_{dR}\alpha|_{\ker \alpha}\otimes_A \mathrm{id}_{H^0(A)}: \ker \alpha \otimes_A H^0(A) \rightarrow (\ker \alpha)^{\vee}[k] \otimes_A H^0(A). $

We first observe that, at $p\in \spec H^0(A)$, $d_{dR}\alpha|_{p}$ maps \begin{align*}
	&\langle \partial/\partial x_1^{-i}|_{p}, \dots, \partial/\partial x_{m_i}^{-i}|_{p} \rangle_{\K} \xrightarrow{\sim} \langle d_{dR}y_1^{k+i}|_{p},\dots, d_{dR}y_{m_i}^{k+i}|_{p} \rangle_{\K}, \\ &\langle \partial/\partial y^{k+i}_1|_{p},\dots,\partial/\partial y^{k+i}_{m_i}|_{p} \rangle_{\K} \xrightarrow{\sim} \langle d_{dR}x_1^{-i}|_{p},\dots, d_{dR}x_{m_i}^{-i}|_{p} \rangle_{\K}
\end{align*} isomorphically for all $i$, such that $d_{dR} \alpha (\partial/\partial z^k|_{p},-)=0$. It follows that, at $p\in \spec H^0(A)$, we have the identifications
\begin{equation}\label{identification T_A/ker d_DR alpha}
	\Big(\faktor{\mathbb{T}_A}{\ker d_{dR}\alpha}\Big)|_p  \simeq \langle \partial/\partial x_j^{-i}|_p, \partial/\partial y^{k+i}_j|_p : \forall i,j\rangle_{\K}		\simeq (\ker \alpha)|_p.
\end{equation}  Thus, the maps $(d_{dR}\alpha|_{\ker \alpha})^{i}|_p: (\ker \alpha|_p)^{i}\rightarrow (\ker^{\vee} \alpha|_p)^{k+i} $ are all isomorphisms at $p$, hence isomorphisms in a neighborhood of $p$. So, localizing $A$ at $p$ if necessary,  the induced morphism $d_{dR}\alpha|_{\ker \alpha}\otimes_A \mathrm{id}_{H^0(A)}$ is an isomorphism of complexes, and hence a quasi-isomorphism. Thus, $d_{dR}\alpha |_{\ker \alpha}$ is non-degenerate, and we get  a \bfem{$k$-shifted contact structure} on $\spec A$ with a $k$-contact form $\alpha$ in the sense of Definition \ref{defn_shiftedcontact} as desired. 
	  
	\begin{definition}
	  	If  $A, d, \alpha$ are as above, 	we then say $A, \alpha$  are \bfem{in contact Darboux form}. 
	  \end{definition}
  Note that the general  expressions like $``d_{dR} z^k + \phi/k +\dR [\cdots]/k"$   will still be valid for the other cases $ (a) \ k\equiv0 \mod 4, \text{ and } (b) \ k\equiv2 \mod 4  $ as well. In fact, Equations (\ref{new local variables for k=-4l}) $-$ (\ref{defn_phiv2}) in Appendix \ref{app. non-odd shifts} show that the other cases involve modified versions of  $H, d, $ and  $\phi$ with some possible extra terms. In either case, the modified $A, \alpha$ would also serve as the desired contact model. Following the same terminology as above, we would also say \bfem{$A, \alpha$ are in (contact) Darboux form}.
\end{example}
\begin{observation}
	Sanity check: the cases $k=-1 \ (\ell=0)$ and $k=-3 \ (\ell=1)$. 
	
\begin{itemize}
	\item 	\bfem{When $k=-1,$} we set $A=A(0)[z^{-1}, y_1^{-1}, \dots, y_{m_0}^{-1}]$, with $A(0)$ a $\mathbb{K}$-algebra generated by $x_1^0,\dots, x_{m_0}^0$, such that $\vdim A= m_0 - (m_0+1)=-1$. Choosing an arbitrary Hamiltonian $H\in A(0)$, we let $dz^{-1}=H$ and $dy^{-1}_j= {\partial H}/{\partial x_j^{0}}$ and $dx^0_j=0 \ \forall j$.  Then, from (\ref{the local primitive contact model}),  the element  $$\alpha=\dR z^{-1}+\sum_{1\leq j \leq m_0}y_j^{-1}\dR x_j^{0}$$ defines a $(-1)$-contact form. In that case, $\ker \alpha$ is generated by the vector fields $ \partial/\partial y_j^{-1}$ { and } $y_j^{-1}\partial/\partial z^{-1} -   \partial/\partial x^{0}_j  \text{ for } 1\leq j\leq m_0.  $ 
	
%	Note that if we let $dz^{-1}=H \in A(0)$, then $d\alpha=0$ implies  the equations $dy^{-1}_j= {\partial H}/{\partial x_j^{0}}$ and $dx^0_j={\partial H}/{\partial y_j^{-1}}=0$ (as $H$ is independent of $y^{-1}_j$'s), which determine $d$ uniquely. Conversely, when $d$ is given by these equations, one has $d\alpha=0$. % as $H$ does not involve any of $y^{-1}_j$'s due to degree reasons.
	
	\item \bfem{When $k=-3,$} we let $A=A(0)[x_1^{-1}, \dots, x_{m_1}^{-1}, y_1^{-2}, \dots, y_{m_1}^{-2}, z^{-3}, y_1^{-3}, \dots, y_{m_0}^{-3}]$, with $A(0)$ a $\mathbb{K}$-algebra generated by $x_1^0,\dots, x_{m_0}^0$, such that $\vdim A=-1$. 
	
	Observe that from \cite[Example 5.17]{Brav}, a Hamiltonian $H\in A^{-2}$ is an element of the form $H=\sum_{j=1}^{m_1}y_j^{-2}s_j+ \sum_{i,j=1}^{m_1}x^{-1}_ix^{-1}_jt_{ij}$, with $s_j,t_{ij}\in A(0)$, satisfying $\sum_{i,j=1}^{m_1}t_{ij}s_j=0$ in $A(0)$ for $i=1,\dots, m_1.$ 
	
	From (\ref{defn_internal d contact}), the differential $d$ is  defined  by  $dx^0_j=0$,  $ dx^{-1}_j=s_j $,   $ dy^{-2}_j=2\sum_{j'=1}^{m_1}x^{-1}_{j'}t_{jj'} $, \ $dy^{-3}_{j}=\sum_{j'=1}^{m_1}y^{-2}_{j'}\frac{\partial s_{j'}}{\partial x^0_{j}} + \sum_{j'',j'} x^{-1}_{j''}x^{-1}_{j'}\frac{\partial t_{j''j'}}{\partial x^0_{j}}$, and \ $ dz^{-3}=\frac{H}{3}-\frac{1}{3} \sum_js_jy^{-1}_j + \frac{2}{3}\sum_{j,j'} x^{-1}_jx^{-1}_{j'}t_{jj'}$.
	Then  the element   $$\alpha=\dR {z}^{-3} +\sum_{j=1}^{m_0} {y}_j^{-3} \dR {x}_j^{0} + \sum_{j=1}^{m_1} {y}_j^{-2}\dR {x}_j^{-1}$$ is a $ (-3) $-contact form. Here, $\ker \alpha$ is  generated by the vector fields $\ y^{-3}_j\partial/\partial z^{-3} -\partial/\partial x^{0}_j$  \text{ and } $\partial/\partial y^{-3}_j$ for $1\leq j\leq m_0;$ \text{and by }  $\partial/\partial x_j^{-1}-y_j^{-2}\partial/\partial z^{-3} \text{ and }  \partial/\partial y_j^{-2}$  for all $1\leq j\leq m_1$.
\end{itemize}
\end{observation}

\begin{example} \label{example_alternatve contact forms}
	Let $k, A, H,\phi$ be as in Example \ref{model example}. Alternatively, we can define the \bfem{differential} $d$ on $A$ as in (\ref{defn_internal d contact}), but with $-kdz^k=H$, instead. In that case, we then introduce the element $ \alpha' \in \Omega_A^1[k]$  by \begin{align*}
		\alpha'&= d_{dR}z^k+  \sum_{i,j} \Big[-\dfrac{i}{k} \ x_j^{-i}d_{dR}y_j^{k+i}+\frac{k+i}{k} \ y_j^{k+i}d_{dR}x_j^{-i}\Big] \\ &=d_{dR}z^k+\phi/k.
	\end{align*}Note that  $d\alpha'=  0$ as well due to the new choice of $dz^k$. Modifying Example \ref{model example} accordingly, one can also conclude that such element $\alpha'$ also serves as a $k$-contact form, inducing a $k$-contact structure on $\spec A$, with the contact data constructed similarly. Details are left to the reader. Using the same terminology, we also say \bfem{$A, \alpha'$ are in (contact) Darboux form}. 
\end{example}
Notice that  both versions  $\alpha$ and $\alpha'$ (and the corresponding differentials) in Example \ref{model example} and Example \ref{example_alternatve contact forms} coincide for $k=-1$. These examples also suggest  that suitable modifications using $H+d[\cdots]$ and $\phi+\dR [\cdots]$ may  lead to alternative versions of such forms. %Note also that using suitable $H+d[\cdots]$ and $\phi+\dR [\cdots]$, one can actually provide further alternative versions of $\alpha$. %in Example \ref{model example}. In any case, we use the same terminology and refer to them as \bfem{(contact) Darboux forms}.
%Following the same terminology, we refer to  all possible such models as \bfem{standard contact Darboux forms}.
\subsubsection{Main results for negatively shifted contact derived schemes} \label{section_darboux theorem}
We now outline the main results of \cite{kib} for  derived $\K$-schemes (of locally finite presentations) with negatively shifted contact structures. 
\paragraph{A Darboux-type theorem.} The first result is about Darboux-type local models, which essentially says that for $k<0$, every $k$-shifted contact derived $\K$-scheme $\bf X$ is locally equivalent to $(\spec A, \alpha_0)$ for $A$ a minimal standard form cdga and $\alpha_0$ as in Example \ref{model example}. More precisely, we have:

\begin{theorem} \cite[Thm. 3.13.]{kib}\label{contact darboux}
	Let $\bf X$ be  a $k$-shifted contact derived $\mathbb{K}$-scheme  for $k<0$, and $x\in \bf X$. Then there is a local contact model $\big(A,  \alpha_0 \big)$  and $p \in \spec H^0(A)$ such that $i: \spec A \hookrightarrow \bf X$ is an open inclusion with $i(p)=x$, \ $A$ is a standard form that is minimal at  $p$, and $\alpha_0$ is a $k$-shifted contact form on $\spec A$ such that $A, \alpha_0$ are in Darboux form.%, with $i^*(\alpha) \sim \alpha_0$ in the space of $k$-shifted 1-forms. 

\end{theorem} Note that for $k<0$ odd, for instance, the pair $(A, \alpha_0)$ can be explicitly given by local graded variables as in Examples \ref{model example} and \ref{example_alternatve contact forms}. For the other cases, one should use another sets of variables as in Equations (\ref{new local variables for k=-4l}) and (\ref{new local variables for k=-4l-2}), and modify $H, \phi, d$ accordingly.
%\paragraph{Symplectization of a shifted contact derived scheme.} \label{section_symplectization}
\paragraph{Symplectification.} The second main result of \cite{kib} is about the  \textit{symplectification} of a $k$-shifted contact derived $\K$-scheme. Recall from \cite[Def. 4.3]{kib} that if $\bf X$ is a locally finitely presented derived $\K$-scheme carrying a $k$-shifted contact structure $(\mathcal{K}, \kappa, L)$ with $k<0$, then we define its \emph{symplectification} $ \mathcal{S}_{{\bf{X}}} $ to be the total space $\widetilde{\bf L}$ of the $\mathbb{G}_m$-bundle of $L$ over $\bf X$, provided with a canonical $k$-shifted symplectic structure (for which the $\mathbb{G}_m$-action is of weight 1) as defined below. 
\newpage
Let $({\bf X}; \mathcal{K}, \kappa, L )$ be a $k$-shifted contact derived $\K$-scheme of locally finite presentation. Given $k<0$ and $p\in \bf X$, find an affine derived sub-scheme ${\bf U}:= \spec A$ such that $p: \spec A \rightarrow \bf X$ is Zariski open inclusion (we may further assume $A$ is of minimal standard form). Here, we assume w.l.o.g. that  $L$ is trivial on $\bf U$.  Define the functor $ \mathcal{S}_{{\bf{X}}}: cdga_{\K} \rightarrow Spcs$ by $ A\mapsto \mathcal{S}_{{\bf{X}}}(A), $ where
\begin{equation}\label{defn_presstack S_X(A)}
	\mathcal{S}_{{\bf{X}}}(A):= \big\{(p, \alpha, v ) :  p \in {\bf X} (A), \ \alpha:   p^*(\mathbb{T}_{{\bf X}})\rightarrow \mathcal{O}[k], \ v: Cocone(\alpha) \xrightarrow{\sim}p^*(\mathcal{K}) \big\},
\end{equation} where each $ v $ is a quasi-isomorphism respecting the natural morphisms $ p^*\kappa: p^*\mathcal{K} \rightarrow p^*(\mathbb{T}_{{\bf X}})$ and $ Cocone(\alpha) \rightarrow p^*(\mathbb{T}_{{\bf X}})   $.
Under the current assumptions, the perfect complexes $\mathbb{T}_A, \mathbb{L}_A$, when restricted to ${\spec H^0(A)}$, are both quasi-isomorphic to free  complexes of $H^0(A)$-modules. For $A \in cdga_{\K}$, we then define a $\mathbb{G}_m(A)$-action on $ \mathcal{S}_{{\bf{X}}}(A)  $ by 
\begin{equation*}
	f \triangleleft (p, \alpha, v):= (p, f\cdot \alpha,v).
\end{equation*}

By \cite[Prop. 4.4]{kib}, $ \mathcal{S}_{{\bf{X}}}$ is equivalent to the total space $\bf \widetilde{L}$ of the $\mathbb{G}_m$-bundle of $L$. Therefore, it has the structure of a derived stack together with the projection map $\pi_1:  \mathcal{S}_{{\bf{X}}}\rightarrow \bf X$.  

We  also introduce the \emph{canonical 1-from} $\lambda$ on $\mathcal{S}_{{\bf X}}$. By construction, we have the projection maps $\pi_1:  \mathcal{S}_{{\bf{X}}}\rightarrow \bf X$ and $\pi_2:  \mathcal{S}_{{\bf{X}}}\rightarrow {\bf T^*[k]X}.$ 
We  define the \emph{canonical 1-from} $\lambda$ on $\mathcal{S}_{{\bf X}}$ to be the pullback  $ \pi_2^*\lambda_{\bf X} $ of the tautological 1-form $ \lambda_{\bf X} $ on $ {\bf T^*[k]X}$.  We then have:

\begin{theorem} \cite[Thm. 4.7]{kib} \label{thm_Symplectization}
	Let $\bf X$ be a (locally finitely presented) derived $\K$-scheme carrying a $k$-shifted contact structure $(\mathcal{K}, \kappa, L)$ with $k<0$.
	Then the $k$-shifted closed 2-form $  {\omega}:= (d_{dR}\lambda, 0, 0, \dots)$ is non-degenerate, and hence the derived stack $ \mathcal{S}_{{\bf{X}}}$ is $k$-shifted symplectic. %which is canonically determined by the shifted contact structure of $\bf X$ (up to quasi-isomorphism). 
	
		We then call the pair $(\mathcal{S}_{\bf{X}}, {\omega})$  the \emph{symplectification of $ \bf X. $}
\end{theorem}

\section{Results for derived Artin stacks} \label{section_Artin stacs}

In this section, we will  explain how to extend the main results of \cite{kib}, outlined in the previous section (cf. Theorems \ref{contact darboux} and \ref{thm_Symplectization}), from derived schemes to the more general case of derived Artin stacks. We begin with some basic definitions and results from \cite{BenBassat}. Later, we give two theorems about the desired generalizations (cf. Theorems \ref{cor_darboux for Artin} \& \ref{cor_symplectization for Artin}). 
\subsection{A Darboux-type theorem for shifted symplectic derived Artin stacks} \label{section_PTVV on Artin}

\paragraph{Nice atlases for derived Artin stacks.} Recall that derived schemes has nice local models. With the same spirit, derived Artin $\K$-stacks  have \textit{nice atlases}. In that respect, we have the following generalization of Definition \ref{defn_standard forms}:
\begin{definition}\label{defn_standard form open nbhd}
	Let $\bf X$ be a derived Artin $\K$-stack and $x\in \bf X$. By a \bfem{standard form open neighborhood of $x$}, we mean a pair $(A, \varphi)$ and a point $p\in \spec H^0(A)$ such that $A$ is a standard form cdga in the sense of Definition \ref{defn_standard forms}, and  $\varphi: {\bf U}=\spec A \rightarrow \bf X$ is  smooth of some relative dimension $n \geq 0$ with $\varphi(p)=x.$
	\end{definition}

For  $A, \ \varphi: {\bf U}=\spec A \rightarrow \bf X$, and a point $p\in \cspec H^0(A)$ as above, there exists a canonical distinguished triangle
\begin{equation} \label{triangle}
\varphi^* \mathbb{L}_{{\bf X}} \rightarrow \mathbb{L}_{{\bf U}} \rightarrow \mathbb{L}_{{\bf U}/{\bf X}} \rightarrow \varphi^* \mathbb{L}_{{\bf X}} \ [1].
\end{equation} As $\varphi$ is  smooth of some relative dimension $n \geq 0$, $ \mathbb{L}_{{\bf U}/{\bf X}} $ is locally free of $ \rank n$. Moreover, for $\varphi(p)=x$, an element of the pullback $\varphi^* \mathbb{L}_{\bf X}$ is locally of the form $ (f \otimes \beta)|_{p}$ with $\beta \in \mathbb{L}_X|_{x}$ and $f\in A$, such that the map $ \varphi^* \mathbb{L}_{{\bf X}}|_p \rightarrow \mathbb{L}_{{\bf U}} |_p $ sends $ f \otimes \beta \mapsto f\cdot\beta.$ %the pullback on stalks gives \begin{equation}
%\varphi^* \mathbb{L}_{{\bf X}} |_p \simeq \mathbb{L}_X |_x \otimes A.
%\end{equation}

\begin{observation}
It follows from the sequence (\ref{triangle}) that $ H^i(\mathbb{L}_{{\bf X}}|_{x}) \simeq H^i(\mathbb{L}_{{\bf U}} |_p) $ for $i<0.$ Moreover, as $\bf U$ is not ``stacky", it follows that $H^1(\mathbb{L}_{\bf{U}}|_p)=0.$ Hence, there exists an exact sequence of $\K$-vector spaces 
\begin{equation} \label{triangle2}
0 \rightarrow H^0(\mathbb{L}_{{\bf X}} |_x) \rightarrow H^0(\mathbb{L}_{{\bf U}} |_p) \rightarrow H^0(\mathbb{L}_{{\bf U}/{\bf X}} |_p) \rightarrow H^1( \mathbb{L}_{{\bf X}} |_x) \rightarrow 0,
\end{equation}so that $n \geq \dim H^1( \mathbb{L}_{{\bf X}} |_x)$ due to the exactness.
\end{observation}
\begin{definition}\label{defn_minimality for std form open nbhd}
		We say that a standard form open neighborhood $(A, \varphi, p)$ of $x$ is \bfem{minimal} if $ A $  is minimal in the sense of Definition \ref{1st defn of minimalty}, and the relative dimension $n$ attains its minimum; i.e., $n=\dim H^1(\mathbb{L}_{\bf X}|_x)$.
\end{definition}

\begin{observation} \label{observation_isomorphisms from triangle}
	Given a minimal standard form open neighborhood $(A, \varphi, p)$ of $x$, %$\varphi$ is smooth of minimal possible dimension $n=\dim H^1(\mathbb{L}_{\bf X}|_x)$, and by 
	the sequence  (\ref{triangle2}) implies that there are isomorphisms $  H^0(\mathbb{L}_{{\bf U}/{\bf X}} |_p) \simeq H^1(\mathbb{L}_{{\bf X}}|_{x})$ and $ H^0(\mathbb{L}_{{\bf X}}|_{x}) \simeq H^0(\mathbb{L}_{{\bf U}} |_p). $ Therefore, we have \begin{equation} \label{isomorphism from triangle}
H^i(\mathbb{L}_{{\bf X}}|_{x}) \simeq H^i(\mathbb{L}_{{\bf U}} |_p) \text{ for } i\leq 0.
	\end{equation} It follows that $A(0)$ is smooth of dimension $m_0= \dim H^0(\mathbb{L}_{\bf{X}}|_{x})$, and $A$ is free over $A(0)$ with $m_i=\dim H^{-i}(\mathbb{L}_{{\bf X}}|_{x})$ generators in each degree $-i$.
\end{observation}

Ben-Bassat, Brav, Bussi and Joyce \cite{BenBassat} proved that derived Artin $\K$-stacks have nice local models in terms of minimal standard form open neighborhoods. The following result summarizes key observations from \cite[Theorems 2.8 \& 2.9]{BenBassat} and in fact serves as a generalization of Theorem \ref{localmodelthm}. For more details, we refer to \cite[Sections 2.4 \& 2.5]{BenBassat}.

\begin{theorem}  \label{thm_localmodelforArtincase}
Let $\bf X$ be a derived Artin $\K$-stack and $x\in \bf X$. Then there exists  a minimal standard form open neighborhood $(A, \varphi, p)$ of $x$.
Moreover, if $(A, \varphi, p)$ and $(A', \varphi', p')$ are two such open neighborhoods, then there exists another standard form cdga $ A'' $ which can be used to compare them in a reasonable way.
\end{theorem}

\paragraph{Darboux form atlases for negatively shifted symplectic derived Artin stacks.}For  $k<0$, it has been  proven in \cite{BenBassat} that given a $k$-shifted symplectic derived Artin $\K$-stack $(\bf X, \omega)$, near each $x\in \bf X$, one can find a ``minimal smooth atlas" $\varphi: \bf U \rightarrow \bf X$ with ${\bf U}=\spec A$ an affine derived scheme such that $({\bf U}, \varphi^*(\omega))$ is in a standard Darboux form. More precisely, we have:

\begin{theorem} (\cite[Theorem 2.10]{BenBassat}) \label{thm_darboux for Artin}
Let $(\bf X, \omega)$ be a $k$-shifted symplectic derived Artin $\K$-stack for $k<0$, and $x\in \bf X$. Then there exist a minimal standard form open neighborhood $(A, \varphi, p)$ of $x$, with $\varphi(p)=x$,  a minimal standard form cdga $B$ with  inclusion $\iota: B \hookrightarrow A$ and the diagram 
\begin{equation} \label{diagram with atlas}
\spec B={\bf V}\xleftarrow{j:=\spec (\iota)} {\bf U}=\spec A \xrightarrow{\varphi} \bf X
\end{equation} such that the induced morphism $\tau({\bf U})\xrightarrow{\tau(j)} \tau ({\bf V})$ between truncations is an isomorphism, and there is a $k$-shifted symplectic structure $\omega_B=(\omega_B^0, 0, 0, \dots)$ on ${\bf V}=\spec B$, which is in Darboux form in the sense of Theorem \ref{Symplectic darboux}, with $\varphi^*(\omega)\sim j^*(\omega_B)$ in $k$-shifted closed 2-forms on $\bf U.$

Moreover, there exists a natural equivalence 
\begin{equation}
\mathbb{L}_{\bf {U/V}} \simeq \mathbb{T}_{\bf{U/V}} [1-k].
\end{equation}
\end{theorem}

\subsection{A Darboux-type theorem for shifted contact derived Artin stacks} \label{section_Contact on Artin}

In this section, we will discuss how to extend  Theorem \ref{contact darboux}  from  derived schemes to derived Artin stacks. In that respect, the following result provides a Darboux-type atlas for \emph{negatively} shifted contact derived Artin stacks. The proof will be a variation of \cite[Theorem 2.10.]{BenBassat}.

\begin{theorem} \label{cor_darboux for Artin}
Given $k\in \mathbb{Z}_{<0}$, let  $\bf X$ be a  derived Artin $\K$-stack (locally of finite presentation) carrying a $k$-shifted contact structure, and $x\in \bf X$. Then we can find \begin{itemize}
	\item  a minimal standard form open neighborhood $(A, \varphi: {\bf U} \rightarrow {\bf X}, p)$ of $x$ \ (cf. Defn.'s \ref{defn_standard form open nbhd} \& \ref{defn_minimality for std form open nbhd});
	\item  a dg-subalgebra $B$ of $A$ with inclusion $\iota: B \hookrightarrow A$ and the diagram 
	\begin{equation*}
		{\bf V}:=\spec B\xleftarrow{j:=\spec (\iota)} {\bf U}=\spec A \xrightarrow{\varphi} \bf X; \ \text{ and }
	\end{equation*}
\item  a $k$-contact form $\alpha_B$ on $ {\bf V}$ such that for any $k$-contact form $\alpha$ on $\bf X$, we have an equivalence $\varphi^*(\alpha)\sim j^*(\alpha_B)$ in $\mathcal{A}^1({\bf U}, k)$, and the pair $({B},\alpha_B)$ is in  contact Darboux form  (cf. Example \ref{model example} for $k$ odd).  %with $(except in degree $k-1$) {\bf U}=\spec A$ an affine derived scheme equipped with a canonical contact data, such that $({\bf U}, \varphi^*(\alpha))$ is in a contact Darboux form for any $k$-contact form $\alpha$ on $\bf X$.
\end{itemize} Moreover,  the induced morphism $\tau(j): \tau({\bf U})\rightarrow{ } \tau ({\bf V})$ between truncations is an isomorphism.
\end{theorem}
\pf Let $(\mathcal{K} \xrightarrow{\kappa} \mathbb{T}_{{\bf X}}, L)$ be a $k$-shifted contact data on $\bf X.$
Then, locally on $\bf X$, where $L$ is trivialized, the perfect complex $\mathcal{K}$ can be  given as a cocone of $ \alpha$, where   $\alpha: \mathbb{T}_{{\bf X}} \rightarrow \mathcal{O}_{{\bf X}}[k]$ is a $k$-contact form, such that we have the (co)exact triangle $\mathcal{K} \rightarrow \mathbb{T}_{{\bf X}} \rightarrow L[k]$. %Throughout the proof, we will use "refined"  neighborhoods introduced in $\S \ref{section_Artin stacs}$; and once the local data is specified, we will fix the corresponding locally defining 1-form.
 
Given $k<0$ and $x\in \bf X$, apply now Theorem \ref{thm_localmodelforArtincase} to get a minimal standard form open neighborhood $(A, \varphi: {\bf U} \rightarrow {\bf X}, p)$ of $x$, with ${\bf U}=\spec A, \ x \in {\bf X}(A), \ p\in \spec H^0(A)$ such that $p\mapsto x$ %Here, from Yoneda's lemma, ${\bf X}(A) \simeq Map_{dPstk}(\spec A, {\bf X})$, and hence any $A$-point $x\in {\bf X}(A) $ can be seen as a morphism $\varphi: \spec A \rightarrow \bf X$ of derived pre-stacks. 
and that the map $\varphi$ is smooth of relative dimension $n=\dim H^1(\mathbb{L}_{\bf X}|_x)$. 
 %and thus its pullback map $ \varphi^*: QCoh({\bf X}) \rightarrow QCoh(\spec A) \simeq Mod_A$ sends $E \mapsto p^*E$. Hence, $s$ is just an element of the $A$-module $ p^*E,$ a "fiber" over $p$.
We also assume that $A$ is a standard form cdga constructed inductively as described in (\ref{A(n) construction})  such that $\mathbb{L}_{A}$ has Tor-amplitude in $[k-1,0]$.

\paragraph{When the shift is odd.} As before,  we again focus on a particular and the simplest case: $k$ is \textit{odd}, say $k=-2\ell-1$ for $\ell\in \N.$  By Definition \ref{defn_standard form open nbhd}, $A$ can be chosen as a free algebra over $A(0)$ with $m_i$ generators in degree $-i$ for $i=1, \dots, \ell$, and $m_i$ generators in degree $k+i$ for $i=0, \dots, \ell$ , but  with additional $n$ generators in degree $k-1$. See Example \ref{example_Artin}.

% and $A(0)$ is a smooth $\mathbb{K}$-algebra of $m_0=\dim H^0(\mathbb{L}_A|_p)$. Moreover, $A$ will have $m_i=\dim H^{-i}(\mathbb{L}_A|_p)$ generators in each degree $-i.$
W.l.o.g., we can assume\footnote{Otherwise, apply again Theorem \ref{thm_localmodelforArtincase} to the overlap $\varphi({\bf U}) \times_{{{\bf X}}}^h W$, where $W$ is open containing $x$ s.t. $L|_W$ is trivial, to get another minimal standard form open neighborhood $(A', \varphi': {\bf U'} \rightarrow {\bf X}, p')$. Then we consider  the interior of $ \varphi'({\bf U'})$ over which the restriction of $L$ is still trivial. Here $ \varphi({\bf U}) $ denotes the image $\Ima \varphi$, with the monomorphism $\Ima \varphi \hookrightarrow \bf X.$} 
that  $L$ is trivial on the interior of the image $\varphi(\bf U)$, denoted by $ \widetilde{\bf U}$. Then, over $ \widetilde{\bf U} $, the induced 1-form $\alpha: \mathbb{T}_{{\bf X}} \rightarrow \mathcal{O}_{{\bf X}}[k]$ is such that $\mathcal{K}$ is the cocone of $\alpha$, up to quasi-isomorphism, and the 2-form $d_{dR}\alpha$ is non-degenerate on $\mathcal{K}$. In that case, we have the exact triangle $\mathcal{K} \rightarrow \mathbb{T}_{{\bf X}} \rightarrow L[k]$  over $ \widetilde{\bf U} $. We fix this  $k$-contact form $\alpha$ for the rest of the proof. Also, by abuse of notation, we simply use $\alpha$ for its pullback on $\mathcal{K}$ as well.

Next, we consider $ d_{dR}\varphi^*(\alpha)$ as a sequence $(d_{dR}\varphi^*(\alpha), 0, 0, \dots),$ which defines a closed $k$-shifted 2-form on ${\bf U}=\spec A$. Applying Proposition \ref{Proposition_exactness} to $d_{dR}\varphi^*(\alpha)$, we obtain  elements $H\in A^{k+1}$ and $\phi \in (\Omega^1_{A})^k$ such that $dH=0$, \ $d_{dR}H+d\phi=0$ , and $kd_{dR}\varphi^*(\alpha) \sim (d_{dR}\phi, 0, 0, \dots).$ We denote this representative by $ (\omega^0, 0, 0, \dots) $ or just $\omega^0$ whenever the meaning is clear.

Note that the pullback $\varphi^*(\alpha)$ may not be a $k$-shifted contact form on $\bf U$, because $ d_{dR}\varphi^*(\alpha) $  is \textit{not} necessarily non-degenerate on $\varphi^*(\mathcal{K})$. But,  we may ensure the non-degeneracy for some degrees. In this regard, we have the following lemma.

\begin{lemma} \label{lemma_q-isom only for some i}
 Let $\omega^0 \sim d_{dR}\varphi^*(\alpha)$ be the  representative of $ d_{dR}\varphi^*(\alpha)$ as above. Then the induced morphism $ \omega^0|_{\varphi^*(\mathcal{K})} \cdot: \varphi^*(\mathcal{K}) \rightarrow \varphi^*(\mathcal{K}^{\vee}[k])$ is a quasi-isomorphism \emph{only for} $0\leq i \leq -k$ \ (i.e.  except in $\deg \ k-1$).

\end{lemma}
\pf[Proof of Lemma \ref{lemma_q-isom only for some i}]
 We first note that by Observation \ref{observation_isomorphisms from triangle}, we have $ H^i(\mathbb{L}_{{\bf X}}|_{x}) \simeq H^i(\mathbb{L}_{\bf{U}} |_p) \text{ for } i\leq 0,$ where $  \mathbb{L}_{{\bf U}}\simeq  \mathbb{L}_{{A}}\simeq\Omega^1_A$ with $ A $ a (minimal) standard form cdga. Then the natural morphism $ \mathbb{L}_{\varphi}[k]: \varphi^*(\mathbb{L}_{{\bf X}})[k]\rightarrow \mathbb{L}_{A}[k] $ from the triangle (\ref{triangle}) induces an isomorphism on cohomology $H^i$ at $p$ for $i+k\leq0$, and it is zero if $i+k=1$. Likewise, the dual $ \mathbb{L}_{\varphi}^{\vee} $ induces an isomorphism on cohomology $H^i$ at $p$ for $i\geq0$, and it is zero if $i=-1$. Thus, we get the  diagram\footnote{The diagonal map $\gamma$ is the composition defined by the commuting triangle.}
 \begin{equation}\label{diagram_nondegeneracy for pullback}
 \begin{tikzcd}
 	T_A \arrow[r, "\mathbb{L}_{\varphi}^{\vee}"]               & \varphi^*(\mathbb{T}_{{\bf X}}) \arrow[rr, "\varphi^*(d_{dR}\alpha \ \cdot)"]                                                            &  & {\varphi^*(\mathbb{L}_{{\bf X}}[k])} \arrow[d, "{\varphi^*\kappa^{\vee}[k]}"] \arrow[r, "{\mathbb{L}_{\varphi}[k]}"] & {\Omega^1_A [k]} \arrow[d, "{\gamma^{\vee}[k]}", dashed] \\
 	{\varphi^*(\mathcal{K}[-1])} \arrow[u] \arrow[r, "\simeq"] & \varphi^*(\mathcal{K}) \arrow[rr, "\omega^0|_{\varphi^*(\mathcal{K})} \ \cdot",dashed] \arrow[u, "\varphi^*\kappa"'] \arrow[lu, "\gamma"', dashed] &  & {\varphi^*(\mathcal{K}^{\vee}[k])} \arrow[r, equal]                                                     & {\varphi^*(\mathcal{K}^{\vee}[k]).}                 
 \end{tikzcd}
 \end{equation}Combining the conditions for $i$ above, we observe that  for $0\leq i \leq -k$ \emph{only}, Diagram \ref{diagram_nondegeneracy for pullback} describes a suitable factorization of the map $\omega^0\cdot: \mathbb{T}_A \rightarrow \Omega^1_A[k]$ commuting 
\begin{equation}
	\begin{tikzcd}
		T_A\simeq \varphi^*(\mathbb{T}_{{\bf X}}) \arrow[rr, "\omega^0 \cdot"] &  & {\Omega^1_A [k]\simeq \varphi^*(\mathbb{L}_{{\bf X}}[k])} \arrow[d] \\
		\varphi^*(\mathcal{K}) \arrow[rr, "\omega^0|_{\varphi^*(\mathcal{K})} \cdot"] \arrow[u]   &  & {\varphi^*(\mathcal{K}^{\vee}[k]),}                         
	\end{tikzcd}
\end{equation}where the 2nd row  is the pullback of $ {\mathcal{K}} \rightarrow {\mathcal{K}^{\vee}}[k] $, which is  an equivalence as  $d_{dR}\alpha$ is non-degenerate on $\mathcal{K}$. Thus, we conclude that the map $\omega^0\cdot: \mathbb{T}_A \rightarrow \Omega^1_A[k]$ representing $d_{dR}\varphi^*(\alpha) \cdot$ is non-degenerate on $ \varphi^*(\mathcal{K}) $ only for $0\leq i \leq -k$.

\epf

 \bfem{Variables in degrees $0, -1, \dots, k$.} Now, we start with the pullback under $\varphi$ of the triangle $\mathcal{K} \rightarrow \mathbb{T}_{{\bf X}} \rightarrow L[k]$ above, which locally splits over $\bf U$. That is, $\mathrm{cofib}(\varphi^*(\mathcal{K}) \rightarrow \varphi^*(\mathbb{T}_{{\bf X}}))\simeq \varphi^*(L[k])$, which is concentrated in $\deg -k$, and $ \varphi^*(\mathbb{T}_{{\bf X}})=  \varphi^*(\mathcal{K})  \oplus \varphi^*(L[k]).$ Here, we will call the 2nd summand $Rest$. Localizing $A$ at $p$ if necessary, first choose  degree 0 variables ${x}_1^0, x_2^0, \dots, x_{m_0}^0$ in $A(0)$ such that $\{d_{dR}x_j^0: j=1,\dots,m_0 \}$ forms a basis for  $\varphi^*(\mathcal{K}^{\vee})^0$ over $A(0),$ and $(Rest^{\vee})^0$ is zero.  
 
 Since the induced morphism given by $d_{dR}\varphi^*(\alpha)$ is an equivalence {on} $ \varphi^*(\mathcal{K}) $ \emph{only  for} $k\leq -i\leq 0$ due to Lemma \ref{lemma_q-isom only for some i}, we can make the following choices of  variables in degrees $0, -1, \dots, k$:  \begin{itemize}
 	\item When $i=0$, we find generators $ y_1^{k}, y_2^{k}, \dots, y_{m_{0}}^{k}, z^{k}  \in A^{k}$ such that $ \{ d_{dR}y_1^{k}, \dots, d_{dR}y_{m_0}^{k}\}$ is a basis for $\varphi^*(\mathcal{K}^{\vee})^{k}$  which is dual to  the basis $\{ d_{dR} x_1^{0}, \dots ,d_{dR} x_{m_0}^{0}\}$ for $\varphi^*(\mathcal{K}^{\vee})^{0}$, and that the complex $Rest$ is  generated by the vector field $\partial/\partial z^{k}$ of degree $-k$.  %Moreover, $d_{dR}\varphi^*(\alpha)\cdot$ maps $ \partial/\partial z^{k} \mapsto 0. $
 	\item For $1 \leq i\leq \ell$, again by  Lemma \ref{lemma_q-isom only for some i}, we can choose generators (except in $\deg \ k-1$) $ x_1^{-i}, x_2^{-i}, \dots, x_{m_i}^{-i}  \in A^{-i}$ and $ y_1^{k+i}, y_2^{k+i}, \dots, y_{m_i}^{k+i}  \in A^{k+i}$ such that $d_{dR}y_j^{k+i}, \text{ for } j= 1, \dots, m_i,$ form a basis for $\varphi^*(\mathcal{K}^{\vee})^{k+i}$ which is dual to  the basis $\{ d_{dR} x_1^{-i}, \dots ,d_{dR} x_{m_i}^{-i}\}$ for $\varphi^*(\mathcal{K}^{\vee})^{-i}$. 
 \end{itemize} %$\dR y_1^{k}\in Rest^{\vee}[k]$\footnote{It is concentrated in $\deg k$ and generated by $ \dR y_1^{k}. $ Note that $\partial/\partial y_1^k$ is a vector field of degree $0$, and hence we must have $\partial/\partial y_1^k\in \varphi^*(\mathcal{K})^0$ by construction.} %is the image of $\partial/\partial \tx_1^{0}$ under the induced map $d_{dR}\varphi^*(\alpha)\cdot: \mathbb{T}_A \rightarrow \mathbb{L}_A[k]$, 
Then the map $\omega^0|_{\varphi^*(\mathcal{K})} \cdot \otimes \ \mathrm{id}_{H^0(A)}$ gives an equivalence\footnote{Due to the minimality at $p$, $d^{-i}|_p=0=(d^{-i})^{\vee}|_p$ for each $i$. So, all degree-wise maps are isomorphisms at $p$, and hence in a neighborhood of $p$. So, localizing $A$ at $p$ if necessary, we can assume $\omega^0|_{\varphi^*(\mathcal{K})} \cdot \otimes \ \mathrm{id}_{H^0(A)}$ is an isomorphism.} of $H^0(A)$-modules in each degree:
\begin{align}
\langle \partial/\partial x_1^{-i}, \dots, \partial/\partial x_{m_i}^{-i} \rangle  &\xrightarrow{\sim} \langle d_{dR} y^{k+i}_1, \dots, d_{dR} y^{k+i}_{m_i} \rangle & &\text{for } 0\leq i\leq \ell, \label{sending duals 1}\\
\langle \partial/\partial y^{k+i}_1, \dots, \partial/\partial y^{k+i}_{m_i} \rangle  &\xrightarrow{\sim} \langle d_{dR} x_1^{-i}, \dots, d_{dR}x^{-i}_{m_i} \rangle & &\text{for } 0\leq i\leq \ell. \label{sending duals 2} 
\end{align} 

Note that we have $\varphi^*(\mathcal{K})=\varphi^*(Cocone(\alpha))\simeq Cocone(\varphi^*(\alpha))$ using the  equivalence between two  exact sequences 
$ Cocone(\varphi^*\alpha)\rightarrow \varphi^*\mathbb{T}_{{\bf X}} \xrightarrow{\varphi^*\alpha} \mathcal{O}[k] \ \text{ and } 
\varphi^*Cocone(\alpha)\rightarrow \varphi^*\mathbb{T}_{{\bf X}} \xrightarrow{\varphi^*\alpha} \mathcal{O}[k], $
where the latter is the pullback of $ Cocone(\alpha)\rightarrow \mathbb{T}_{{\bf X}} \xrightarrow{\alpha} \mathcal{O}_{{\bf X}}[k]. $ Moreover, over $\spec A$, we  simply write $ \ker \varphi^*(\alpha) $ instead of $Cocone(\varphi^*\alpha)$\footnote{Thanks to Proposition \ref{defn_contact strc on good affines}.}.

Thus, using the local coordinates, the splitting $ \mathbb{T}_A=  \ker \varphi^*(\alpha) \oplus Rest$  is such that $ \ker \varphi^*(\alpha) $ has $\mathrm{Tor}$-amplitude $[0,-k]$ and $Rest$ is concentrated in $\deg -k$,  where
\begin{align} \label{explicit generators for ker and rest}
 \ker \varphi^*(\alpha) |_{\spec H^0(A)}&= \big\langle  \partial/\partial x_j^{-i}, \partial/\partial y^{k+i}_j: 0\leq i\leq \ell, \ 1\leq j \leq m_i \big\rangle_{H^0(A)},  \nonumber \\
Rest|_{\spec H^0(A)} &=\big\langle \partial/\partial z^{k} \big\rangle_{H^0(A)}. 
\end{align}
Then using the complexes in (\ref{explicit generators for ker and rest}), the non-degeneracy condition for $d_{dR}\varphi^*(\alpha)$ on $ \ker \varphi^*(\alpha) $ \emph{- except for $\deg \  (k-1)$ -} sending the dual basis of $ d_{dR} x_a^{b}$ to the basis $  d_{dR} y^{b'}_{a'}$ (and vice versa) as in Equations (\ref{sending duals 1}) - (\ref{sending duals 2}) implies that $ d_{dR}\varphi^*(\alpha)\in \wedge^2\Omega_A^1[k] $ is given by
\begin{equation}\label{2 form over kernel}
	d_{dR}\varphi^*(\alpha)=\displaystyle \sum_{i=0}^{\ell} \sum_{j=1}^{m_i} d_{dR}x_j^{-i} d_{dR}y_j^{k+i}.
\end{equation}
Note  that the kernel  of the  induced map $d_{dR}\varphi^*(\alpha) \cdot$ is  spanned by the degree $-k$ vector field $\partial/\partial z^k$, while the action  of $d_{dR}\varphi^*(\alpha) \cdot$ on $ \ker \varphi^*(\alpha) $ is given by Equations (\ref{sending duals 1})-(\ref{sending duals 2}). Thus, we  obtain the identification $\faktor{\mathbb{T}_{{A}}}{\ker \dR\varphi^*(\alpha)}\simeq  \ker \varphi^*(\alpha)$.

Scaling $z^k$ we may assume $\iota_{\partial/\partial z^k} \varphi^*(\alpha) =1$. Now, our goal is to find a unique\footnote{The conditions uniquely determine the explicit form, up to interchange of $x_{j}^{-i} \text{ and } y_{j}^{k+i}$. Here, the roles of $x_{j}^{-i},y_{j}^{k+i}$ are symmetric in (\ref{2 form over kernel}) and (\ref{explicit generators for ker and rest}), where $ \dR x_{j}^{-i}\dR y_{j}^{k+i}= \dR y_{j}^{k+i}\dR x_{j}^{-i} $ for $k$ odd. See \cite[Proof of Thm. 3.13]{kib}.} $\varphi^*(\alpha)$ satisfying Eqn. (\ref{2 form over kernel}), the condition on the kernel in (\ref{explicit generators for ker and rest}), and the equation $\iota_{\partial/\partial z^k} \varphi^*(\alpha) =1 $. %Since the  form $ \varphi^*(\alpha) $ can also been seen as a representative of $\phi$, we may let $ \varphi^*(\alpha)=\dR z^k + \phi/k $. 
 Then such  $ \varphi^*(\alpha) $ satisfying the desired properties can be  written explicitly as
\begin{equation} \label{desired form}
 \varphi^*(\alpha)= d_{dR}z^k+  \sum_{i,j} y_j^{k+i}d_{dR}x_j^{-i}.
\end{equation}  %Observe that $ d_{dR}\varphi^*(\alpha) $ recovers (\ref{2 form over kernel}) and that the kernel of $ \varphi^*(\alpha) $ for the above representation is generated by the the vectors fields\footnote{The vector fields $ ix^{-i}_j\partial/\partial x_j^{-i} + (k+i)y^{k+i}_j\partial/\partial y^{k+i}_j $ also annihilate $\varphi^*(\alpha)$, but they can be expressed using the others.} $ k \partial/\partial x_j^{-i}-(k+i)y^{k+i}_j\partial/\partial z^{k} \text{ and } k\partial/\partial y^{k+i}_j+ix^{-i}_j\partial/\partial z^{k} $, which can be identified with the ones in (\ref{explicit generators for ker and rest}).  

%\vspace{2pt}

\bfem{Variables in degree $(k-1)$ and the differential.} We construct the rest by combining Example \ref{model example} with Example \ref{example_Artin}. Note that  Example \ref{model example} does \emph{not} involve the additional finitely many generators in degree $(k-1)$. However, due to the atlas chosen at the beginning of the proof, the corresponding cdga $A$ must admit additional $n$ generators in degree $(k-1)$   as in Example \ref{example_Artin}.

In our case, we  identify $A$ as a \emph{commutative graded algebra} with the commutative graded algebra  freely generated by the variables $z^{k}, x_j^{-i}, y_{j'}^{k+i'}$ as above, but \bfem{ with additional $n$ generators, $ w_1^{k-1}, \dots, w_{n}^{k-1}$, in degree $k-1$}. As discussed before, $\omega^0, H, \phi$ above do not involve any of $ w_j^{k-1} $ for degree reasons, and hence the extra variables can be chosen arbitrarily. %Likewise, $ dw_j^{k-1}$ can be  determined arbitrarily. %We can also impose a suitable differential $d$ as before: $d$ acts on $z^{k}, x_j^{-i}, y_{j'}^{k+i'}$ as in Equation (\ref{defn_internal d}); and $ dw_j^{k-1}$ can be  determined arbitrarily.

Choose $B$  with  inclusion $\iota: B \hookrightarrow A$   such that  $B(0)$ is the subalgebra of $A(0)$ with the same generators $ x_1^0, \dots, x_{m_0}^0$ and that the sub-cdga $B$ is  the free algebra over $B(0)$ on the generators $x_{j}^{-i}, y_j^{k+i},z^k$ \bfem{only}. That is, we identify $B$ as a \emph{commutative graded algebra} with the commutative graded algebra in Example \ref{model example}.

 It remains to show that  $H$ satisfies (\ref{defn_CME}) and the \bfem{differential $d$ on $B$} can be given by Equation (\ref{defn_internal d contact}). To this end, we analyze the defining equations for the pair $(H,\phi)$. Notice that $d$ will not be fully determined  on $A$, but determined only on $B$, which is enough for our construction. So, $ dw_j^{k-1}$ can be arbitrary. 

First of all,   combining the  defining equation  $d_{dR}\phi=kd_{dR}\varphi^*(\alpha)$ above with the equation (\ref{2 form over kernel}),  we may explicitly write\footnote{See the proof of \cite[Theorem 5.18]{Brav}. Alternatively,  one can let $ \phi=  k\sum_{i=0}^{\ell} \sum_{j=1}^{m_i} y_j^{k+i}d_{dR}x_j^{-i}$.  Leaving $\dR\phi$ unchanged, these expressions can be transformed to each other by replacing $H, \phi$ by suitable $H+d(\cdots), \ \phi+d_{dR}(\cdots),$ respectively. See Examples \ref{model example} and \ref{example_alternatve contact forms}.}  $$ \phi= \sum_{i=0}^{\ell} \sum_{j=1}^{m_i} [(-i)x_j^{-i}d_{dR}y_j^{k+i}+(k+i)y_j^{k+i}d_{dR}x_j^{-i}].$$ Then (the proof of ) \cite[Theorem 5.18]{Brav}  shows that  expanding 
$ d_{dR}H+d\phi=0 $ with the explicit representation of $\phi$ above and comparing the coefficients of $d_{dR}$-terms, one can get the following formulas for $d$:
\begin{equation} 
	d|_{B(0)}=0; \ dx_j^{-i} =  \dfrac{\partial H}{\partial y_j^{k+i}} \text{ for all } i>0,j; \ \text{ and } \ dy_j^{k+i} =  \dfrac{\partial H}{\partial x_j^{-i}} \text{ for all } i,j.
\end{equation} 

Secondly, using these equations\footnote{These equations will be enough as $H$ is independent of the variables $z^k, y^k_j$'s.}  to expand $dH=0$, \cite[Theorem 5.18]{Brav}  also shows that $H$ satisfies the classical master equation (\ref{defn_CME}). 
Before the final step, we also observe that using the explicit representation of $\phi$ above, $\varphi^*(\alpha)$  in (\ref{desired form}) can also be rewritten as %(cf. Example \ref{model example})
  \begin{equation} \label{form written alternatively}
	\varphi^*(\alpha)=\dR z^k + \dfrac{1}{k}\Big[\phi + \dR  \big[\sum_{i,j} (-1)^{i} ix_j^{-i} y_j^{k+i} \big]\Big].
\end{equation}

Finally, combining the defining equation $ d_{dR}H+d\phi=0 $ with Eqn. (\ref{form written alternatively}) (and $d\varphi^*(\alpha)=0$), we get $ \dR H=-d\phi=-kd\varphi^*(\alpha) + kd\circ\dR z^k +d\circ\dR [\cdots]= \dR (-kdz^k - d[\cdots]).  $ So, we  have such $ H $ satisfying $ -kdz^k=H+d[\cdots]$. 
Thus, we conclude that the differential $d$ is given by Eqn. (\ref{defn_internal d contact});  hence, $(B,d)$ is identified with the cdga in Example  \ref{model example}.

\vspace{3pt}

\bfem{Contact data on $\spec B  \text{ and }  \spec A$.} We  let ${\bf V}:=\spec B$, along with the diagram 
\begin{equation}
	\spec B={\bf V}\xleftarrow{j:=\spec (\iota)} {\bf U}=\spec A \xrightarrow{\varphi} \bf X.
\end{equation} 
For degree reasons, $H\in B$, and $\omega^0, \phi$ are all images under $\iota$ of $\omega_B^0, \phi_B$, respectively, where \begin{equation}
	\omega^0_B=\sum_{i=0}^{\ell} \sum_{j=1}^{m_i} d_{dR}x_j^{-i} d_{dR}y_j^{k+i} \text{ and } \phi_B=\sum_{i=0}^{\ell} \sum_{j=1}^{m_i} [(-i)x_j^{-i}d_{dR}y_j^{k+i}+(k+i)y_j^{k+i}d_{dR}x_j^{-i}].
\end{equation}
Note that from Example \ref{model example}, %$(B,d)$ is a contact  standard form cdga with coordinates $x_{j}^{-i}, y_j^{k+i},z^k$, for $0\leq i \leq \ell; 1 \leq j\leq m_i,$ and with $d$  given by Equation (\ref{defn_internal d contact}), such that 
$\spec B$ carries the canonical $k$-shifted contact structure defined by $\alpha_B=d_{dR}z^k+  \sum_{i,j}y_j^{k+i}d_{dR}x_j^{-i}$.
Therefore,  from Eqn. (\ref{desired form}), we conclude that $$\varphi^*(\alpha)\sim j^*(\alpha_B) \text{ except in degree } (k-1).$$

Lastly, we can consider $j:\bf U \rightarrow V$ as an embedding of $\bf U$ into $\bf V$ as a derived subscheme; and hence,  the induced morphism $\tau(j): \tau({\bf U})\rightarrow{ } \tau ({\bf V})$ between truncations is an isomorphism. This completes the proof for odd $k.$

\paragraph{When the shift is not odd.} For the other cases $ (a) \ k\equiv0 \mod 4, \text{ and } (b) \ k\equiv2 \mod 4  $, one should use another sets of variables as in Equations (\ref{new local variables for k=-4l}) and (\ref{new local variables for k=-4l-2}), respectively, along  with the additional generators $ w_j^{k-1} $ in degree $k-1$. %and modify $H, \phi, d$ as in Equations (\ref{new local variables for k=-4l}) $ - $ (\ref{defn_phiv2})
We leave  details to the reader.

\epf

\subsection{Symplectifications of shifted contact derived Artin stacks}

In this section, %combining the  constructions in the proofs of Theorem \ref{thm_Symplectization} and Theorem \ref{cor_darboux for Artin}, 
we describe the canonical symplectification of  a (negatively) shifted contact derived Artin $\K$-stack. More precisely, we have:

\begin{theorem} \label{cor_symplectization for Artin}
Let  $\bf X$ be a $k$-shifted contact derived Artin $\K$-stack, then its \emph{symplectization} $(\mathcal{S}_{\bf{X}}, \omega)$ can be canonically described as in Theorem \ref{thm_Symplectization}.
\end{theorem}

\pf

First of all, using Theorem \ref{thm_Symplectization}, we define the space $\mathcal{S}_{{\bf X}}$ as follows. Let $({\bf X}; \mathcal{K}, \kappa, L )$ be a $k$-shifted contact derived Artin $\K$-stack of locally finite presentation.  Then we define the functor $ \mathcal{S}_{{\bf{X}}}: cdga_{\K} \rightarrow \mathbb{S}$ by $ A\mapsto \mathcal{S}_{{\bf{X}}}(A), $ where
\begin{equation}
	\mathcal{S}_{{\bf{X}}}(A):= \big\{(x, \alpha, v ) :  x \in {\bf X} (A), \ \alpha:   \varphi^*(\mathbb{T}_{{\bf X}})\rightarrow \mathcal{O}[k], \ v: Cocone(\alpha) \xrightarrow{\sim}\varphi^*(\mathcal{K}) \big\},
\end{equation} where each $ v $ is a quasi-isomorphism respecting the natural morphisms $ \varphi^*\kappa: \varphi^*\mathcal{K} \rightarrow \varphi^*(\mathbb{T}_{{\bf X}})$ and $ Cocone(\alpha) \rightarrow \varphi^*(\mathbb{T}_{{\bf X}})   $.

By \cite[Prop. 4.4]{kib},
$ \mathcal{S}_{{\bf{X}}}$ is equivalent to the total space $\bf \widetilde{L}$ of the $\mathbb{G}_m$-bundle of $L$. Therefore, it has the structure of a derived stack together with the projection maps $\pi_1:  \mathcal{S}_{{\bf{X}}}\rightarrow \bf X$ and $\pi_2:  \mathcal{S}_{{\bf{X}}}\rightarrow {\bf T^*[k]X}.$ 
We  also define the \emph{canonical 1-from} $\lambda$ on $\mathcal{S}_{{\bf X}}$ to be the pullback  $ \pi_2^*\lambda_{\bf X} $ of the tautological 1-form $ \lambda_{\bf X} $ on $ {\bf T^*[k]X}$.

We set $  {\omega}:= (d_{dR}\lambda, 0, 0, \dots)$, which is a $k$-shifted closed 2-form on  $\mathcal{S}_{{\bf X}}$, and hence it defines a pre-$k$-shifted symplectic structure on  $\mathcal{S}_{{\bf X}}$.
 Now, it remains to show that $  {\omega}$ is non-degenerate. As before, the rest of the argument is in fact local, so it is enough to prove it using suitable local models studied in the previous sections. 

Let us now analyze our local data. Given $k<0$ and $x\in \bf X$, find an affine derived scheme ${\bf U}:= \spec A$ and $p \in \spec H^0(A)$ with $\varphi: \spec A \rightarrow \bf X$  a smooth map of relative dimension $n$  such that $p\mapsto x$ (we may further assume $A$ is of minimal standard form due to Theorem \ref{thm_localmodelforArtincase}). Here, we assume w.l.o.g. that  $L$ is trivial on the interior, $\bf \tilde{U}$, of $\bf \Ima \varphi$ as before. Under the current assumptions, the perfect complexes $\mathbb{T}_A, \mathbb{L}_A$, when restricted to ${\spec H^0(A)}$, are both quasi-isomorphic to free  complexes of $H^0(A)$-modules. For $A \in cdga_{\K}$, we then define a $\mathbb{G}_m(A)$-action on $ \mathcal{S}_{{\bf{X}}}(A)  $ as before.

By Theorem \ref{cor_darboux for Artin}, we can construct a cdga $B$  with  inclusion $\iota: B \hookrightarrow A$ and the diagram  \begin{equation*}
	\spec B={\bf V}\xleftarrow{j:=\spec (\iota)} {\bf U}=\spec A \xrightarrow{\varphi} \bf X
\end{equation*}such that the induced morphism $\tau(j): \tau({\bf U})\rightarrow{ } \tau ({\bf V})$ between truncations is an isomorphism.    In fact,  we get the pair $({\bf V}, \alpha_B)$ with an equivalence $\varphi^*(\alpha)\sim j^*(\alpha_B)$ except in degree $k-1$ such that $\alpha_B$ is a $k$-shifted contact form on $\bf V$ as in Example \ref{model example}.

Fixing the local data above, we then consider the homotopy pullback diagram
\begin{equation} 
	\begin{tikzpicture}
		\matrix (m) [matrix of math nodes,row sep=3em,column sep=5 em,minimum width=3 em] {
			\mathcal{Z}:= {\bf U} \times_{\bf X}^h \mathcal{S}_{{\bf{X}}}	& \mathcal{S}_{{\bf{X}}}  & \bf T^*[k]X.  \\			
			{\bf U}  & \bf X   &  \\
		};
		\path[-stealth]
		
		%horizontal 1st row
		(m-1-1) edge  node [above] { $pr_2$} (m-1-2)
		(m-1-2) edge  node [above] { $\pi_2$} (m-1-3)
		
		%horizontal 2nd row
		(m-2-1) edge  node [above] { $ \varphi $ } (m-2-2)
		
		%vertical
		(m-1-1) edge  node [left] {$ pr_1 $ } (m-2-1)
		(m-1-2) edge  node [right] { $ \pi_1 $ } (m-2-2)
	
		%edge [dashed,-] (m-2-1)
		;
	\end{tikzpicture} 
\end{equation} Notice that, over $p\in \bf U$, we can then identify the fiber over $p$ locally as $ {\bf {U}} \times_{{{\bf X}}}^h \mathbb{G}_m,$ with natural projections, because $ \mathcal{S}_{{\bf{X}}}$ is identified the total space of $L$. 

On the part of the space $\mathcal{S}_{\bf{X}}$ over ${\bf \tilde{U}}$, for  the elements $\pi_1^*\alpha, \lambda \in \mathbb{L}_{\mathcal{S}_{\bf{X}}}[k]$, the identification of $\mathcal{S}_{\bf{X}}$ with the total space of $L$ (i.e. the space of trivializations) implies that there is an element $f\in \mathbb{G}_m(A)$ such that  \begin{equation}
\lambda=f \cdot \pi_1^* (\alpha).
\end{equation} %with $f\in \mathbb{G}_m(A).$ %We again assume that $f$ has inverse (or apply localization if necessary). 
%Note also that $\lambda=f \cdot \pi^* \varphi^*(\alpha) \sim f \cdot \pi^*j^*(\alpha_B)$, except in degree $k-1$.
Using the homotopy $\varphi \circ pr_1 \sim \pi_1\circ pr_2$, we get an element $\tilde{\lambda}\in pr_2^*(\mathbb{L}_{\mathcal{S}_{\bf{X}}}[k]) $ such that
\begin{equation}
	\tilde{\lambda}:= pr_2^*(\lambda)=pr_2^*(f)\cdot (\pi_1 \circ pr_2)^*(\alpha) \sim pr_2^*(f)\cdot (\varphi \circ pr_1)^*(\alpha)= pr_2^*(f)\cdot pr_1^*(\varphi^*\alpha),
\end{equation} where we denote $ pr_2^*(f), pr_1^*(\varphi^*\alpha) $ simply by $ \tilde{f}, \widetilde{\varphi^*\alpha}$, respectively. Thus, we get a local representative of $\lambda$, $$\tilde{\lambda}=\tilde{f}\cdot\widetilde{\varphi^*\alpha},$$ on a (minimal) standard form open neighborhood $(A, \varphi, p)$ of $x$.% with an atlas $\varphi: {\bf U}=\spec A\rightarrow {\bf X} \text{ such that } p\mapsto x$.

Moreover, for the map $ \varphi \circ pr_1: \mathcal{Z} \rightarrow \bf X$ we have an exact triangle 
\begin{equation}
(\varphi \circ pr_1)^* (\mathbb{L}_{{\bf X}}[k]) \rightarrow pr_1^*(\mathbb{L}_{{\bf U}}[k]) \oplus pr_2^*(\mathbb{L}_{\mathcal{S}_{\bf{X}}}[k]) \rightarrow \mathbb{L}_{\mathcal{Z}}[k].
\end{equation}
%Given $k<0$ and $x\in \bf X$, apply Theorem \ref{thm_localmodelforArtincase} to get a minimal standard form open neighborhood $(A, \varphi: {\bf U} \rightarrow {\bf X}, p)$ of $x$ with ${\bf U}=\spec A$ and  $p\in \cspec H^0(A)$. 
Now, to prove $  {\omega}:= (d_{dR}\lambda, 0, 0, \dots)$ is non-degenerate, it suffices to show that $d_{dR}\tilde{\lambda}$ is non-degenerate.

\begin{lemma}\label{lemma_non-ged for local reps}
$d_{dR}\tilde{\lambda}$ is non-degenerate.
\end{lemma}
\pf
Recall first that there exists a natural equivalence  $DR(A)\otimes_{\K} DR(C) \simeq DR(A\otimes_{\K} C)$ induced by the identification
\begin{equation} \label{cotangent cmpx of tensor}
	\mathbb{L}_{A\otimes_{\K}C}\simeq(\mathbb{L}_A  \otimes_{\K} C) \oplus (A \otimes_{\K} \mathbb{L}_C). 
\end{equation} In our case, since the fiber over $p$ is locally given as $ {\bf {U}} \times_{{{\bf X}}}^h \mathbb{G}_m,$ where  $ {\bf U}=\spec A $ and $\mathbb{G}_m=\spec C$ is the affine group scheme, with say $C:=\K [x,x^{-1}]$, we can use the identification (\ref{cotangent cmpx of tensor}) to locally decompose $ \tilde{\lambda}=\tilde{f}\cdot\widetilde{\varphi^*\alpha}$.

We then have, when restricted to $\spec H^0(A)$, 
\begin{equation} \label{splitting symplectization}
	(\mathbb{L}_{A\otimes_{\K}C})^{\vee} \simeq \big(\ker (\varphi^*\alpha)\footnote{Recall that on refined local charts, we  simply write $ \ker (\varphi^*\alpha) $ instead of $Cocone(\varphi^*\alpha)$.}  \oplus Rest \big) \oplus \big(H^0(A) \otimes_{\K}\langle \partial/\partial \tilde{f} \rangle_{C}\big).
\end{equation}
Now, to prove that $d_{dR} \tilde{\lambda}$ is  non-degenerate, 
it suffices to show that, at $p\in \bf U,$ for any non-vanishing (homogeneous) vector field $\sigma\in (\mathbb{L}_{A\otimes_{\K}C})^{\vee}$, there is  a vector field $\eta \in (\mathbb{L}_{A\otimes_{\K}C})^{\vee}$ such that $\iota_{\eta} (\iota_{\sigma} d_{dR} \tilde{\lambda}) \neq 0.$\footnote{This argument follows from the following observations: When restricted to $\spec H^0(A)$, the induced morphism $ \mathbb{T}_{A} \rightarrow \Omega^1_{A}[k]$ is just a map of finite complexes of free $H^0(A)$-modules. And, at $p\in \spec H^0(A)$, both $ \mathbb{T}_{A}|_p,  \Omega^1_{A}|_p$ are complexes of $\K$-vector spaces. For non-degeneracy, we require  this map to be a (degree-wise) quasi-isomorphism. Recall that localizing $A$ at $ p $ if necessary, we may assume that the induces map is indeed an (degree-wise) isomorphism near $p$. Therefore, the fact we use here is just an analogous result from linear algebra.} To this end, we first compute
\begin{equation*}
	\iota_{\eta} (\iota_{\sigma} d_{dR} \tilde{\lambda})= \mp (d_{dR}\tilde{f})(\sigma) \varphi^*\alpha(pr_{1,*}\eta) \mp (d_{dR}\tilde{f})(\eta) \varphi^*\alpha(pr_{1,*}\sigma)  \mp \tilde{f}\cdot d_{dR}(\varphi^*\alpha) (pr_{1,*} \sigma, pr_{1,*} \eta).
\end{equation*} From Equation (\ref{splitting symplectization}), it is enough to consider the following cases:
\begin{enumerate}
	\item If $\sigma \in \ker (\varphi^*\alpha)$,  then %$ \pi_* \sigma \simeq \sigma $ and 
	$ \iota_{\eta} (\iota_{\sigma} d_{dR} \tilde{\lambda}) = \mp \tilde{f} d_{dR}(\varphi^*\alpha) (pr_{1,*}\sigma, pr_{1,*} \eta)$. Since $d_{dR}(\varphi^*\alpha) |_{\ker (\varphi^*\alpha)}$ is non-degenerate by the contactness condition on $\varphi^*(\alpha)\sim j^*(\alpha_B)$ except in degree $k-1$ (and $\tilde{f}\neq 0$ as $f\in \mathbb{G}_m(A)$), it is enough to take $\eta$ to be any non-zero vector in $\ker (\varphi^*\alpha)$.
	\item If $\sigma \in Rest,$ then we get %$ \pi_* \sigma \simeq \sigma $ and 
	$ \iota_{\eta} (\iota_{\sigma} d_{dR} \tilde{\lambda}) =\mp (d_{dR}\tilde{f})(\eta) \varphi^*\alpha(\sigma). $ Observe that $\varphi^*\alpha(\sigma) \neq 0 $ since $\sigma \in Rest.$ Thus, it is enough to take $\eta$ to be any non-zero vector in $ H^0(A) \otimes_{\K}\langle \partial/\partial \tilde{f} \rangle_{C} $ so that $(d_{dR}\tilde{f})(\eta) \neq 0.$
	\item If $\sigma \in H^0(A) \otimes_{\K}\langle \partial/\partial \tilde{f} \rangle_{C},$ then $ \iota_{\eta} (\iota_{\sigma} d_{dR} \tilde{\lambda}) = \mp (d_{dR}\tilde{f})(\sigma) \varphi^*\alpha(pr_{1,*}\eta).$ Note that we have $ (d_{dR}\tilde{f})(\sigma)\neq 0,$ so it suffices to take $\eta$ to be any non-zero vector in $ Rest $ so that $ \varphi^*\alpha(pr_{1,*}\eta)\neq 0. $
\end{enumerate} \epf

In total, from Lemma \ref{lemma_non-ged for local reps}, the $k$-shifted 2-form $\omega^0:= d_{dR} \lambda$ is  non-degenerate (except in degree $k-1$), and hence the sequence $\omega:=(\omega^0, 0, 0, \dots)$ defines a $k$-shifted symplectic structure on $\mathcal{S}_{\bf X}.$  
This completes the proof of Theorem \ref{cor_symplectization for Artin}.

\epf
\section{Examples}\label{section_examples}
In this section we present several constructions of derived Artin stacks with shifted contact structure. In brief, the first example $\S\ref{example_shifted jet space}$ generalizes the 1-jet bundles in the classical setup; and the second set of constructions in $\S\ref{examples_contactvia shifted prequ}$ arises from the  notion of shifted geometric  (pre)quantization introduced by Safronov \cite{Safronov2023}.  
\subsection{Shifted 1-jet  stacks} \label{example_shifted jet space}
	Denote by $\mathbb{A}^1$ the \emph{affine line}\footnote{We may also call it the \emph{affine addtive group scheme $ \mathbb{G}_a.$}} as the derived stack corepresented by $\mathbb{K}[z].$ That is, $\mathbb{A}^1= \spec (\mathbb{K}[z])$ is the derived stack \begin{equation}
		\spec (\mathbb{K}[z]): B\in cdga_{\mathbb{K}} \mapsto Hom(\mathbb{K}[z], B). 
	\end{equation} Let $  T^*[n]X $ be the $n$-shifted cotangent stack of the derived Artin stack $ X$ (locally of finite presentation). Consider  the derived stack ${T^*[n]X} \times \mathbb{A}^1[n]  $ given by the pullback diagram \begin{equation} \label{diagram_cartesian spaces as pullbacs}
	\begin{tikzpicture}
		\matrix (m) [matrix of math nodes,row sep=3em,column sep=5 em,minimum width=2 em] {
			J^1[n]X:=T^*[n]X \times \mathbb{A}^1[n]	   & \mathbb{A}^1[n]  \\
			T^*[n]X   &  *, \\};
		\path[-stealth]
		(m-1-1) edge  node [left] { $ pr_1 $} (m-2-1)
		(m-1-1) edge  node [above] { $ pr_2 $} (m-1-2)
		%(m-1-1) edge  node [below] {} node [below] {{\small stacks}} (m-2-2)
		%(m-1-1) edge  node [below] {} node [below] {{\small  higher stacks}} (m-3-2)
		(m-2-1) edge  node  [above] {$  $ } (m-2-2)
		
		(m-1-2) edge  node [right] {$  $ } (m-2-2);
		%edge [dashed,-] (m-2-1);
	\end{tikzpicture}
\end{equation}where $\mathbb{A}^1[n]$ is the \emph{$n$-shifted affine line} corepresented by the polynomial algebra on a variable in cohomological degree $ -n $. Denote this resulting space by $J^1[n]X$ and call it \bfem{$n$-shifted 1-jet stack of $X$,} with the natural projection map as the composition $J^1[n]X\rightarrow T^*[n]X\rightarrow X.$

Recall that \cite{Calaque2016} showed that there is an $n$-shifted 1-form $\lambda$, called the \emph{Liouville one-form}, on the cotangent stack $ T^*[n]X $ such that  the closed 2-form $\omega:=d_{dR}\lambda$ is non-degenerate. Hence, $ T^*[n]X $ is in fact $n$-symplectic. Moreover, from Diagram \ref{diagram_cartesian spaces as pullbacs}, we have the identifications
\begin{equation}\label{identifications by pullback}
	\mathbb{L}_{J^1[n]X} \simeq pr_1^*\mathbb{L}_{T^*[n]X} \oplus pr_2^* \mathbb{L}_{\mathbb{A}^1[n]} \text{ and } \mathbb{T}_{J^1[n]X} \simeq pr_1^*\mathbb{T}_{T^*[n]X} \oplus pr_2^* \mathbb{T}_{\mathbb{A}^1[n]} .
\end{equation}
 Equivalently, we have the exact triangle \begin{equation}\label{triangle giving precontact}
pr_1^*\mathbb{T}_{T^*[n]X} \rightarrow \mathbb{T}_{J^1[n]X} \rightarrow pr_2^* \mathbb{T}_{\mathbb{A}^1[n]}.
 \end{equation} Notice that from the triangle (\ref{triangle giving precontact}), we have a natural morphism $pr_1^*\mathbb{T}_{T^*[n]X} \rightarrow \mathbb{T}_{J^1[n]X}$ with the cofiber $ pr_2^* \mathbb{T}_{\mathbb{A}^1[n]}. $ Since $ \mathbb{T}_{\mathbb{A}^1[n]} $ is an invertible quasi-coherent sheaf, we get $pr_2^* \mathbb{T}_{\mathbb{A}^1[n]} \simeq L [n]$, where $L$ is a line bundle.

Thus, we obtain a natural \emph{pre-$n$-shifted contact data} by using the perfect complex $\mathcal{K}:= pr_1^*\mathbb{T}_{T^*[n]X}$ and the natural map $\kappa: pr_1^*\mathbb{T}_{T^*[n]X} \rightarrow \mathbb{T}_{J^1[n]X}$ with the cofiber $ pr_2^* \mathbb{T}_{\mathbb{A}^1[n]}$ as above.

Now, it remains to extend the pre-contact data to an $n$-shifted contact structure. In brief, we need to represent $\mathcal{K}$ (locally) as the cocone of a shifted 1-form satisfying the  non-degeneracy condition.  To this end, we construct the following $n$-shifted 1-form on $ J^1[n]X.$

We  define (globally) an $n$-shifted 1-form on $ J^1[n]X $ by\footnote{Thanks to the identifications (\ref{identifications by pullback}).} \begin{equation}
	\alpha:= -d_{dR}z + \lambda \in \Gamma\big(J^1[n]X, \mathbb{L}_{J^1[n]X}[n] \big),
\end{equation} where we simply write $z, \lambda$ instead of $pr_2^*z, pr_1^*\lambda,$ respectively. Then we claim:

\begin{lemma} \label{lemma_useful equivalences for jet stack}
	Let $\alpha, pr_1, pr_2$ be as above. Then we have equivalences (at least locally)
	\begin{align}
		\mathbb{T}_{pr_1}&\simeq pr_2^*\mathbb{T}_{\mathbb{A}^1[n]},\\
		Cocone(\alpha)&\simeq pr_1^*\mathbb{T}_{T^*[n]X}=:\mathcal{K},
	\end{align}where $ \mathbb{T}_{pr_1} $ is the \emph{relative tangent complex}\footnote{Its elements are called \emph{vertical tangent vectors}.} defined by the fiber sequence $\mathbb{T}_{pr_1}\rightarrow \mathbb{T}_{J^1[n]X} \rightarrow pr_1^*\mathbb{T}_{T^*[n]X}$. We then have the (local) splitting \begin{equation*}
\mathbb{T}_{J^1[n]X}\simeq Cocone(\alpha) \oplus \mathbb{T}_{pr_1}.
\end{equation*}
\end{lemma}
\pf
Combining the shift of the natural fiber sequence $\mathbb{T}_{pr_1}\rightarrow \mathbb{T}_{J^1[n]X} \rightarrow pr_1^*\mathbb{T}_{T^*[n]X}$ with the exact triangle (\ref{triangle giving precontact}), we get an equivalence of triangles

\begin{equation}
	\begin{tikzpicture}
		\matrix (m) [matrix of math nodes,row sep=3em,column sep=5 em,minimum width=3 em] {
			pr_1^*\mathbb{T}_{T^*[n]X}	& \mathbb{T}_{J^1[n]X}  & pr_2^*\mathbb{T}_{\mathbb{A}^1[n]}\\
			pr_1^*\mathbb{T}_{T^*[n]X}[-1]	& \mathbb{T}_{J^1[n]X}  & \mathbb{T}_{pr_1}[1],  \\			
		};
		\path[-stealth]
		
		%horizontal 1st row
		(m-1-1) edge  node [above] { $ $} (m-1-2)
		(m-1-2) edge  node [above] { $ $} (m-1-3)
		%(m-1-4) edge  node [right] { } (m-1-5)
		%(m-1-5) edge  node [right] { } (m-1-6)
		%(m-1-6) edge  node [right] { } (m-1-7)
		%(m-1-7) edge  node [right] { } (m-1-8)
		%horizontal 2nd row
		(m-2-1) edge  node [above] { $  $ } (m-2-2)
		(m-2-2) edge  node [above] { $   $ } (m-2-3)
		%(m-2-4) edge  node [right] { } (m-2-5)
		%(m-2-5) edge  node [right] { } (m-2-6)
		%(m-2-6) edge  node [right] { } (m-2-7)
		%	(m-2-7) edge  node [right] { } (m-2-8)
		%vertical
		(m-1-1) edge   node [right] {$ \simeq $ } (m-2-1)
		(m-1-2) edge  node [right] { $ id $ } (m-2-2)
		(m-1-3) edge  [dashed,->] node [right] {$ \simeq $ } (m-2-3)
		%(m-1-6) edge  node [right] { } (m-2-6)
		%(m-1-7) edge  node [right] { } (m-2-7);
		%edge [dashed,-] (m-2-1)
		;
	\end{tikzpicture}
\end{equation}which gives the identification $pr_2^*\mathbb{T}_{\mathbb{A}^1[n]} \xrightarrow{\sim}\mathbb{T}_{pr_1}.$

By the definition of $\alpha$, for any vertical vector $v\in \mathbb{T}_{pr_1}$, the contraction $\iota_v \alpha$ is never nullhomotopic, and hence we write 
$$Cocone(\alpha) \cap \mathbb{T}_{pr_1} = \{0\},$$ by which we mean the pullback of the diagram $Cocone(\alpha)\hookrightarrow \mathbb{T}_{J^1[n]X} \leftarrow \mathbb{T}_{pr_1}$ is trivial in $\mathrm{Perf}(J^1[n]X)$. Since both $Cocone(\alpha), \mathbb{T}_{pr_1}$ are perfect, we then have the (local) splitting\footnote{Over a (minimal) standard form cdga $A$ with an $A$-point $p: \spec A \rightarrow J^1[n]X$, the complex $ p^*\mathbb{T}_{J^1[n]X} $ is a finite complex of free $H^0(A)$-modules, so are the both complexes $Cocone(p^*\alpha), p^*\mathbb{T}_{pr_1}$. Since $ Cocone(p^*\alpha) \cap p^*\mathbb{T}_{pr_1} = \{0\} $ in each degree, we get the splitting of free modules in each degree, and hence the desired (local) identification of complexes.} 
\begin{equation}
	\mathbb{T}_{J^1[n]X} \simeq Cocone(\alpha) \oplus \mathbb{T}_{pr_1}.
\end{equation} 
The splitting then gives an exact triangle
\begin{equation}
Cocone(\alpha) \rightarrow \mathbb{T}_{J^1[n]X} \rightarrow \mathbb{T}_{pr_1},
\end{equation} which induces an equivalence of  triangles
\begin{equation}
	\begin{tikzpicture}
		\matrix (m) [matrix of math nodes,row sep=3em,column sep=5 em,minimum width=3 em] {
			pr_1^*\mathbb{T}_{T^*[n]X}	& \mathbb{T}_{J^1[n]X}  & pr_2^*\mathbb{T}_{\mathbb{A}^1[n]}\\
			Cocone(\alpha)	& \mathbb{T}_{J^1[n]X}  & \mathbb{T}_{pr_1}. \\			
		};
		\path[-stealth]
		
		%horizontal 1st row
		(m-1-1) edge  node [above] { $ $} (m-1-2)
		(m-1-2) edge  node [above] { $ $} (m-1-3)
		%(m-1-4) edge  node [right] { } (m-1-5)
		%(m-1-5) edge  node [right] { } (m-1-6)
		%(m-1-6) edge  node [right] { } (m-1-7)
		%(m-1-7) edge  node [right] { } (m-1-8)
		%horizontal 2nd row
		(m-2-1) edge  node [above] { $  $ } (m-2-2)
		(m-2-2) edge  node [above] { $   $ } (m-2-3)
		%(m-2-4) edge  node [right] { } (m-2-5)
		%(m-2-5) edge  node [right] { } (m-2-6)
		%(m-2-6) edge  node [right] { } (m-2-7)
		%	(m-2-7) edge  node [right] { } (m-2-8)
		%vertical
		(m-1-1) edge   [dashed,->]node [right] {$ \simeq $ } (m-2-1)
		(m-1-2) edge  node [right] { $ id $ } (m-2-2)
		(m-1-3) edge   node [right] {$ \simeq $ } (m-2-3)
		%(m-1-6) edge  node [right] { } (m-2-6)
		%(m-1-7) edge  node [right] { } (m-2-7);
		%edge [dashed,-] (m-2-1)
		;
	\end{tikzpicture}
\end{equation}Thus, we get the identification\footnote{Using the previous footnote, we actually have a local identification $Cocone(p^*\alpha)\xrightarrow{\sim} p^*(pr_1^*\mathbb{T}_{T^*[n]X})= p^*\mathcal{K}.$} $Cocone(\alpha)\xrightarrow{\sim} pr_1^*\mathbb{T}_{T^*[n]X}$ and complete the proof.
\epf
 Finally, using Lemma \ref{lemma_useful equivalences for jet stack}, we can satisfy the desired non-degeneracy condition for $\alpha$ and prove:
 
\begin{theorem}
	Let $ X$ be a locally finitely presented derived Artin stack. Then the $n$-shifted 1-jet stack $ J^1[n]X$  has a natural $n$-shifted contact structure.
\end{theorem}
\pf
Consider the \emph{pre-$n$-shifted contact data} given by the perfect complex $\mathcal{K}:= pr_1^*\mathbb{T}_{T^*[n]X}$ and the natural morphism $\kappa: pr_1^*\mathbb{T}_{T^*[n]X} \rightarrow \mathbb{T}_{J^1[n]X}$ with the cofiber $ pr_2^* \mathbb{T}_{\mathbb{A}^1[n]}$ as above, where  we  have an equivalence  $pr_2^* \mathbb{T}_{\mathbb{A}^1[n]} \simeq L [n]$, with $L$ a line bundle. It remains to check the (local) contact non-degeneracy condition.

Let $\alpha, \lambda$ be as above. By Lemma \ref{lemma_useful equivalences for jet stack}, $ Cocone(\alpha)\simeq pr_1^*\mathbb{T}_{T^*[n]X}=:\mathcal{K}$  locally. Note that $ d_{dR}\lambda $ is already $n$-symplectic on $ T^*[n]X $, and hence non-degenerate on $ \mathbb{T}_{T^*[n]X}.$  Due to the identification by Lemma \ref{lemma_useful equivalences for jet stack}, $d_{dR}\alpha= d_{dR}\lambda$ is then non-degenerate on $\mathcal{K}\simeq Cocone(\alpha)$. Therefore, the data $(\mathcal{K}, \kappa, L, \alpha)$ defines  an $n$-shifted contact structure on $ J^1[n]X$. 

\epf
\subsection{Contact structures via shifted prequantization}\label{examples_contactvia shifted prequ}

\subsubsection{Background on geometric quantization}

\paragraph{Classical geometric quantization.} Let us briefly recall the notion of geometric quantization in the context of differential  geometry. Given a smooth symplectic manifold $ (X, \omega) $ (or a scheme over a field $ \mathbb{K} $ of characteristic 0), \emph{geometric quantization} is a 2-step procedure combining \emph{prequantization} and \emph{polarization.} More precisely, we have: \begin{definition}
	By a \emph{prequantization} of $ (X, \omega) $, we mean a particular line bundle with connection $(L,\nabla)$ on $X$ such that $\mathrm{curv}(L,\nabla)=\omega.$ By a \emph{polarization}, we mean the choice of a subbundle $P$ of $TX$ closed under the Lie bracket which is Lagrangian with respect to $\omega.$
\end{definition}

Note that on $ (X, \omega) $, we can define the  \textit{Hamiltonian vector field $  X_f $ associated to} the function $f$ by $ \imath_{ X_f} \omega=df $, and hence the  \textit{Poisson bracket} 
$ \{f,g \}_{\omega} := -w(X_f,X_g)=X_f (g) $ on  $C^{\infty}(X)$.

When we have a geometric (pre)quantization\footnote{A sufficient condition for the existence is $[\omega]\in H^2(X,\mathbb{Z})$} of $ (X, \omega)$, then it allows us to define  a suitable quantum Hilbert space $\mathcal{H}:=\Gamma_P(X,L)$  as the \emph{space of P-polarizaed sections of $L$} such that we can construct a Lie algebra representation of (a certain subalgebra $A$ of) $(C^{\infty}(X), \{ - , -\}_{\omega})$ on $End(\mathcal{H}, [-,-])$, which acts on the functions as \begin{equation*}
	f\mapsto \hat{f}:=-i\hbar \nabla_{X_f} - f.
\end{equation*}
\begin{example}
	 Every K\"{a}hler manifold $ (M, \omega, J) $ gives rise to \textit{a holomorphic K\"{a}hler polarization} associated to $(M, \omega)$ by setting ${P}:= T^{(0,1)}(M)$, the $(-i)$-eigenspace subbundle of the complexified tangent bundle $TM \otimes \mathbb{C}$.
\end{example}

\begin{remark} \emph{(Relation with contact geometry)}
	The construction above can be given in terms of the principal $U(1)$-bundle $L^{\times}$ associated with $L$ and the connection 1-form $\alpha$ on $X$ corresponding to $\nabla$. We then have $\mathrm{curv}(L^{\times},\alpha)= d_{dR}\alpha=\pi^*\omega$, with $\pi:L^{\times}\rightarrow X$. In that case, $\alpha$ defines a contact structure on $L^{\times}$. For more details, see \cite{Weinstein2005}.
\end{remark}

\paragraph{Shifted geometric quantization.}

Safronov introduces in \cite{Safronov2023} the notion of shifted geometric quantization in the context of derived symplectic geometry. We now outline the main definitions and some key results of interest. We closely follow \cite{Safronov2023}.

Denote by $\mathcal{A}^p(n), \mathcal{A}^{p,cl}(n)$ the \emph{derived stacks of $ p $-forms of degree $ n $} and \emph{closed $ p $-forms of degree $ n $}
as introduced in PTVV's work \cite{PTVV}. Then by construction, we have equivalences \begin{equation} \label{equivalences_loop space vs form}
\mathcal{A}^1(n)\simeq \Omega \mathcal{A}^1(n+1) \text{ and }  	\mathcal{A}^{1,cl}(n)\simeq \Omega \mathcal{A}^{1,cl}(n+1).
\end{equation} By adjunction we can get a map $B\mathcal{A}^1(n)\rightarrow  \mathcal{A}^1(n+1)$, which is an equivalence \cite[Lemma 2.1]{Safronov2023}.

 \begin{observation} \label{observ_some equivalences regarding shifted forms}
 	Iterating these equivalences, we also get an equivalence $$B^n\mathcal{A}^1(0)\rightarrow  \mathcal{A}^1(n).$$ However, it should be noted that $B^n\mathcal{A}^{1,cl}(0)\rightarrow  \mathcal{A}^{1,cl}(n)$ is not an equivalence even if the inclusion $\mathcal{A}^{n+1,cl}(0)\rightarrow  \mathcal{A}^{1,cl}(n)$ is an equivalence. For details, we refer to \cite[$\S$2.1]{Safronov2023}.
 \end{observation}

Note that that we have a morphism of  stacks
\begin{equation}
	d_{dR} \log : \mathbb{G}_m \rightarrow \mathcal{A}^{1,cl}(0)
\end{equation}which, on a cdga $ R $, maps  an invertible element $ f\in R^{\times} $ to $ \dfrac{d_{dR}f}{f} \in \mathcal{A}^{1,cl}(\spec R, 0). $ Then by delooping (cf. Observation \ref{observ_some equivalences regarding shifted forms}), we get the map
\begin{equation}
	c_1:B^n\mathbb{G}_m \rightarrow B^n\mathcal{A}^{1,cl}(0) \rightarrow  \mathcal{A}^{1,cl}(n).
\end{equation}Now, pre-composing $c_1$ with the natural projection $ \mathcal{A}^{1,cl}(n)\rightarrow\mathcal{A}^{1}(n),$ we get a morphism $$  B^n\mathbb{G}_m \rightarrow \mathcal{A}^{1}(n),$$ which will be denoted again by $c_1$ and lead to the following definition.

\begin{definition}
	Let $n\geq 0$ and $X$ a derived stack.
	\begin{enumerate}
		\item An \emph{\textbf{$n$-gerbe}} on $X$ is a map $\mathcal{G}: X \rightarrow B^{n+1}\mathbb{G}_m.$ We then call its   image $c_1(\mathcal{G}) \in \mathcal{A}^{1}(n+1)$ under the map $c_1$ above the \emph{characteristic class of $\mathcal{G}$}
		\item An \emph{\textbf{$n$-gerbe with a connective structure}} on $X$ is a pair $(\mathcal{G}, \nabla)$, where $\mathcal{G}$ is an $n$-gerbe on $X$ and $\nabla$ is a nullhomotopy of the characteristic class $c_1(\mathcal{G}) \in \mathcal{A}^{1}(n+1)$.
		\item An \emph{\textbf{$n$-gerbe with a flat connective structure}} on $X$ is a pair $(\mathcal{G}, \nabla)$, where $\mathcal{G}$ is an $n$-gerbe on $X$ and $\nabla$ is a nullhomotopy of $c_1(\mathcal{G}) \in \mathcal{A}^{1,cl}(n+1)$.
		\item Any gerbe $(\mathcal{G}, \nabla)$ with connective structure on $X$ has a \emph{\textbf{curvature}} $\mathrm{curv}(\mathcal{G}, \nabla) \in \mathcal{A}^{2,cl}(n)$ defined by the pullback diagram
		 \begin{equation}
			\begin{tikzpicture}
				\matrix (m) [matrix of math nodes,row sep=3em,column sep=5 em,minimum width=2 em] {
					 & \mathcal{A}^{2,cl}(n)  & 0 \\
					B^{n+1}\mathbb{G}_m	& \mathcal{A}^{1,cl}(n+1)  & \mathcal{A}^{1}(n+1). \\			
				};
				\path[-stealth]
				
				%horizontal 1st row
				%(m-1-1) edge  node [above] { $ i$} (m-1-2)
				(m-1-2) edge  node [above] { $ $} (m-1-3)
				%(m-1-4) edge  node [right] { } (m-1-5)
				%(m-1-5) edge  node [right] { } (m-1-6)
				%(m-1-6) edge  node [right] { } (m-1-7)
				%(m-1-7) edge  node [right] { } (m-1-8)
				%horizontal 2nd row
				(m-2-1) edge  node [above] { $ c_1 $ } (m-2-2)
				(m-2-2) edge  node [above] { $   $ } (m-2-3)
				%(m-2-4) edge  node [right] { } (m-2-5)
				%(m-2-5) edge  node [right] { } (m-2-6)
				%(m-2-6) edge  node [right] { } (m-2-7)
				%	(m-2-7) edge  node [right] { } (m-2-8)
				%vertical
				%(m-1-1) edge   node [right] {$ id$ } (m-2-1)
				(m-1-2) edge  node [right] { $ $ } (m-2-2)
				(m-1-3) edge  node [right] {$  $ } (m-2-3)
				%(m-1-6) edge  node [right] { } (m-2-6)
				%(m-1-7) edge  node [right] { } (m-2-7);
				%edge [dashed,-] (m-2-1)
				;
			\end{tikzpicture} 
		\end{equation}
	\end{enumerate} 
\end{definition}

\begin{example}
	Let $(\mathcal{G}, \nabla)$ be an $n$-gerbe with connective structure on $X$. When $n=0$, the map $\mathcal{G}: X \rightarrow B^{}\mathbb{G}_m$ corresponds to a line bundle as $B\mathbb{G}_m$ is the classifying space of the line bundles. Then the map $c_1:B\mathbb{G}_m \rightarrow \mathcal{A}^{1,cl}(1)\simeq\mathcal{A}^{2,cl}(0)$ corresponds the first Chern class of the line bundle. Thus, the connective structure $ (\mathcal{G},\nabla)$ corresponds to the usual notion of a connection on the line bundle $\mathcal{G}$ and its curvature corresponds to the usual notion of the curvature. For more examples and details, see \cite{Safronov2023}.
\end{example}

Now, we define the notion of prequantization for shifted symplectic derived stacks.

\begin{definition}
	Let $(X,\omega)$ be an $n$-shifted symplectic derived stack. Its \emph{\textbf{prequantization}} is an $n$-gerbe with connective structure $(\mathcal{G}, \nabla)$ such that $\mathrm{curv}(\mathcal{G}, \nabla) \simeq \omega$ in $ \mathcal{A}^{2,cl}(n). $
\end{definition}

We can also define a relative version of the prequantization notion for the maps $\pi: X\rightarrow B$ of  derived Artin stacks carrying a shifted Lagrangian fibration structure.

\begin{definition}
	Let $\pi: X\rightarrow B$ be a morphism of derived Artin stacks locally of finite presentation. A \emph{\textbf{prequantum $ n $-shifted Lagrangian fibration}} consists of the following data:
	\begin{enumerate}
		\item An $ n $-gerbe $\mathcal{G}$ on $ B $.
		\item An extension of the natural relative flat connection on $\pi^*\mathcal{G}$ to a connective structure $\nabla$, such that  the pair $ (\pi^*\mathcal{G}, \nabla) $ defines an $n$-shifted Lagrangian fibration structure\footnote{For the precise definition, see \cite{Calaque2016}. For our purposes, it is useful to know some key results: Once we have such structure, then $(1)$ the source $X$ is an $n$-symplectic stack; (2) For  any $b\in B$, the inclusion of the fiber $X_b\rightarrow X$ has an $ n $-shifted Lagrangian structure.} on $\pi$ with an induced $n$-shifted symplectic structure on $X$ given by $ \mathrm{curv}(\pi^*\mathcal{G}, \nabla) $.
	\end{enumerate}
\end{definition}
It is also possible to define the prequantization notion for the morphisms $\pi: X\rightarrow B$ of  derived Artin stacks carrying a shifted Lagrangian structure. For details, see \cite[$\S$2.3]{Safronov2023}.

In the upcoming sections, we discuss several examples of prequantizations from \cite{Safronov2023} and provide sample constructions of the induced shifted contact structures.
\subsubsection{Prequantization of the cotangent stack}\label{example_contact induced by prequantum on the cotangent}
Let $X$ be a derived Artin stack locally of finite presentation and $\pi_X: T^*X\rightarrow X$ the natural projection. \cite{Calaque2016} shows that the cotangent stack $T^*X$ has a natural Liouville 1-form $\lambda \in \mathcal{A}^1(T^*X, 0)$ such that 0-shifted symplectic structure on $T^*X$ is given by $\omega:=d_{dR}\lambda$ and that the map $\pi_X$ carries a natural structure of an 0-shifted Lagrangian fibration. That is, for any $x\in X$, the inclusion of the fiber $X_x\rightarrow T^*X$ has a 0-shifted Lagrangian structure.

Now, we consider another interesting structure that $\pi_X$ carries: It is shown in \cite[Prop. 2.21]{Safronov2023} that $\pi_X$ has a natural structure of a prequantum 0-shifted Lagrangian fibration, determined by the trivial 0-gerbe $\mathcal{G}$ on $X$ together with a
connective structure on $\pi_X^*\mathcal{G}$ given by the Liouville form $\lambda$. 

By definition, the trivial 0-gerbe is the map $\mathcal{G}:X\rightarrow B\mathbb{G}_m$ which in fact corresponds to a trivial line bundle, say ${L}$, on $X$. Moreover, in that case, the connective structure $\pi_X^*\mathcal{G}$ corresponds to the usual connection $\nabla$ on the trivial line bundle $\pi_X^*L$ over $T^*X$, with $\pi: \pi_X^*L \rightarrow T^*X,$  such that the connection 1-form is given by $\nabla:= \lambda \in \mathcal{A}^1(T^*X, 0).$

Denote the trivial $\mathbb{G}_m$-bundle associated with $\pi_X^*L$ by $\mathcal{L}^{\times}$. Here, the frame bundle $\mathcal{L}^{\times}$ has the trivialization as the restriction of the trivialization of the original line bundle $\pi: \pi_X^*L \rightarrow T^*X,$ respecting the $\mathbb{G}_m$-action. Likewise, the projection $\pi: \mathcal{L}^{\times} \rightarrow T^*X$ is the restriction of  the original projection $\pi: \pi_X^*L \rightarrow T^*X$.

Since $\mathcal{L}^{\times}$  is a trivial $\mathbb{G}_m$-bundle, we can view it as the derived stack given by the  pullback 

\begin{equation} \label{diagram_G_m bundle as pullbac}
	\begin{tikzpicture}
		\matrix (m) [matrix of math nodes,row sep=3em,column sep=5 em,minimum width=2 em] {
			\mathcal{L}^{\times}=T^*X \times \mathbb{G}_m	   & \mathbb{G}_m  \\
			T^*X   &  *, \\};
		\path[-stealth]
		(m-1-1) edge  node [left] { $ \pi=:pr_1 $} (m-2-1)
		(m-1-1) edge  node [above] { $ pr_2 $} (m-1-2)
		%(m-1-1) edge  node [below] {} node [below] {{\small stacks}} (m-2-2)
		%(m-1-1) edge  node [below] {} node [below] {{\small  higher stacks}} (m-3-2)
		(m-2-1) edge  node  [above] {$  $ } (m-2-2)
		
		(m-1-2) edge  node [right] {$  $ } (m-2-2);
		%edge [dashed,-] (m-2-1);
	\end{tikzpicture}
\end{equation} where $\mathbb{G}_m$ is the affine group scheme as the derived stack corepresented by $\mathbb{K}[t,t^{-1}]$: \begin{equation}
	\mathbb{G}_m=\spec \mathbb{K}[t,t^{-1}] : R \in cdga_{\mathbb{K}} \mapsto Hom(\mathbb{K}[t,t^{-1}], R),
\end{equation} which also means $\mathbb{G}_m(R)=R^{\times}.$ Moreover, from Diagram \ref{diagram_G_m bundle as pullbac}, we have the identifications
\begin{equation}\label{identifications G_m by pullback}
\mathbb{L}_{\mathcal{L}^{\times}} \simeq pr_1^*\mathbb{L}_{T^*X} \oplus pr_2^* \mathbb{L}_{\mathbb{G}_m} \text{ and } \mathbb{T}_{\mathcal{L}^{\times}} \simeq pr_1^*\mathbb{T}_{T^*X} \oplus pr_2^* \mathbb{T}_{\mathbb{G}_m}.
\end{equation}

Recall also from the previous section that we have a morphism of  stacks
\begin{equation}
	d_{dR} \log : \mathbb{G}_m \rightarrow \mathcal{A}^{1,cl}(0)
\end{equation}which, on a cdga $ R $, maps  an invertible element $ f\in R^{\times} $ to $ {(d_{dR}f)}/{f} \in \mathcal{A}^{1,cl}(\spec R, 0). $

For $R=\mathbb{K}[t,t^{-1}]= \mathcal{O}(\mathbb{G}_m)$ and $f:=t$, we get a global section of $\mathcal{L}^{\times}$ as $t \mapsto (d_{dR}t)/t,$ and hence  a \emph{global 0-shifted 1-form} $\alpha$ on $ \mathcal{L}^{\times} $ defined by 
\begin{equation} \label{defn_global 1 form}
\alpha= \pi^* \lambda + pr_2^*d_{dR} \log (t).
\end{equation} By abuse of notation, we may omit $pr^*_1,\pi, pr^*_2$ whenever the meaning is clear from the context.

\begin{observation}
	For the morphism $\pi: \mathcal{L}^{\times} \rightarrow T^*X$, we have the fiber sequence $$\mathbb{T}_{\pi} \rightarrow \mathbb{T}_{\mathcal{L}^{\times}} \rightarrow \pi^* \mathbb{T}_{T^*X},$$ where $\mathbb{T}_{\pi}$ denotes the \emph{relative tangent complex}. An element of the relative tangent space is called a \emph{relative or vertical tangent vector}.
	That is, $ \mathbb{T}_{\pi} $ is  the \emph{space of tangent vectors along the fibers} of $ \mathcal{L}^{\times}.$ Call this space the  \emph{\bf{vertical bundle of} $ \mathbb{T}_{\mathcal{L}^{\times}} $} and denote it by $Ver_{\mathcal{L}^{\times}}.$ Then we have:   
\end{observation}

\begin{lemma} \label{lemma_useful equivalences}
	Let $\alpha, Ver_{\mathcal{L}^{\times}}$ be as above. Then we have the following equivalences: \begin{align}
		Ver_{\mathcal{L}^{\times}}&\simeq \mathrm{cofib} (Cocone(\alpha)\hookrightarrow \mathbb{T}_{\mathcal{L}^{\times}})\simeq pr_2^*\mathbb{T}_{\mathbb{G}_m},\label{lemma_eq1} \\
		 Cocone(\alpha)&\simeq  \pi^*\mathbb{T}_{T^*X}, \label{lemma_eq2}\\
		\mathbb{T}_{\mathcal{L}^{\times}}&\simeq Cocone(\alpha) \oplus  	Ver_{\mathcal{L}^{\times}} \label{lemma_eq3}.
	\end{align} We then call $ Cocone(\alpha) $ the \emph{\textbf{horizontal bundle of} $ \mathbb{T}_{\mathcal{L}^{\times}} $}, denoted by $Hor_{\mathcal{L}^{\times}},$ and write the splitting in terms of the horizontal and vertical bundles as \begin{equation}
\mathbb{T}_{\mathcal{L}^{\times}}\simeq Hor_{\mathcal{L}^{\times}} \oplus  	Ver_{\mathcal{L}^{\times}}.
\end{equation}
\end{lemma}
\pf
Denote the natural morphism $Cocone(\alpha)\hookrightarrow \mathbb{T}_{\mathcal{L}^{\times}}$ by $i$. Then by definition, we have $Cone(i)=\mathrm{cofib} (Cocone(\alpha)\hookrightarrow \mathbb{T}_{\mathcal{L}^{\times}})$, the weak quotient of $ \mathbb{T}_{\mathcal{L}^{\times}} $ by the image of $ Cocone(\alpha).$ 

Notice that for $v\in Ver_{\mathcal{L}^{\times}}$, we have $\pi_*(v) \sim 0$, and hence 
\begin{equation*}
	\alpha (v) = \pi^*\lambda (v) + \frac{1}{t} d_{dR}t (v)= \lambda (\pi_* v) + \frac{1}{t} v (t) \sim \frac{1}{t} v (t).
\end{equation*} Over any $R$-point $p$ of $ \mathcal{L}^{\times} $, $p:\spec R\rightarrow \mathcal{L}^{\times}$, $p^*v$ is a derivation on $R^{\times}$, and hence it maps $t_R:=p^*t$ to an invertible element $p^*v (t_R)$ of $R$. It follows that the restriction of $\alpha$ to the vertical bundle $ Ver_{\mathcal{L}^{\times}} $ will take non-zero values only. I.e., for any $v\in Ver_{\mathcal{L}^{\times}}$, the image $ \alpha (v) $ is homotopic to a non-zero element. Therefore, we can write $$Cocone(\alpha) \cap Ver_{\mathcal{L}^{\times}} = \{0\},$$ by which we mean the pullback of the diagram $Cocone(\alpha)\hookrightarrow \mathbb{T}_{\mathcal{L}^{\times}} \leftarrow \mathbb{T}_{\pi}$ is trivial in $\mathrm{Perf}(\mathcal{L}^{\times})$. Since both $Cocone(\alpha), \mathbb{T}_{\pi}$ are perfect, we then have the (local) splitting\footnote{Over an $A$-point $p$ with $A$ a minimal standard form cdga,  complexes on each side are in fact finite complexes of free $H^0(A)$-modules, and we have the splitting in each degree.} 
\begin{equation}\label{splitting_ker+vertical}
\mathbb{T}_{\mathcal{L}^{\times}} \simeq Cocone(\alpha) \oplus \mathbb{T}_{\pi}.
\end{equation}

By the splitting (\ref{splitting_ker+vertical}), we obtain an equivalence of exact triangles
 \begin{equation}
	\begin{tikzpicture}
		\matrix (m) [matrix of math nodes,row sep=3em,column sep=5 em,minimum width=3 em] {
			Cocone(\alpha) & \mathbb{T}_{\mathcal{L}^{\times}}   & \mathbb{T}_{\pi}\\
			Cocone(\alpha)	& \mathbb{T}_{\mathcal{L}^{\times}}  & Cone(i). \\			
		};
		\path[-stealth]
		
		%horizontal 1st row
		(m-1-1) edge  node [above] { $ i$} (m-1-2)
		(m-1-2) edge  node [above] { $ $} (m-1-3)
		%(m-1-4) edge  node [right] { } (m-1-5)
		%(m-1-5) edge  node [right] { } (m-1-6)
		%(m-1-6) edge  node [right] { } (m-1-7)
		%(m-1-7) edge  node [right] { } (m-1-8)
		%horizontal 2nd row
		(m-2-1) edge  node [above] { $  $ } (m-2-2)
		(m-2-2) edge  node [above] { $   $ } (m-2-3)
		%(m-2-4) edge  node [right] { } (m-2-5)
		%(m-2-5) edge  node [right] { } (m-2-6)
		%(m-2-6) edge  node [right] { } (m-2-7)
		%	(m-2-7) edge  node [right] { } (m-2-8)
		%vertical
		(m-1-1) edge   node [right] {$ id$ } (m-2-1)
		(m-1-2) edge  node [right] { $ id $ } (m-2-2)
		(m-1-3) edge  [dashed,->] node [right] {$ \simeq $ } (m-2-3)
		%(m-1-6) edge  node [right] { } (m-2-6)
		%(m-1-7) edge  node [right] { } (m-2-7);
		%edge [dashed,-] (m-2-1)
		;
	\end{tikzpicture} 
\end{equation}Using the notation $Ver_{\mathcal{L}^{\times}}:=\mathbb{T}_{\pi}$ and the rightmost vertical map on the diagram above, we conclude \begin{equation}
Ver_{\mathcal{L}^{\times}}\simeq Cone(i):=\mathrm{cofib} (Cocone(\alpha)\hookrightarrow \mathbb{T}_{\mathcal{L}^{\times}}),
\end{equation} which proves the first equivalence in (\ref{lemma_eq1}).

Recall that since $\mathcal{L}^{\times}$  is a trivial $\mathbb{G}_m$-bundle, from the identifications (\ref{identifications G_m by pullback}), there is an exact triangle $ 	\pi^*\mathbb{T}_{T^*X}	\rightarrow \mathbb{T}_{\mathcal{L}^{\times}}  \rightarrow pr_2^*\mathbb{T}_{\mathbb{G}_m}$.  Now, using the natural  triangle $\mathbb{T}_{\mathcal{L}^{\times}} \rightarrow \pi^* \mathbb{T}_{T^*X} \rightarrow \mathbb{T}_{\pi}[1]$, we get an equivalence of triangles
\begin{equation}
	\begin{tikzpicture}
		\matrix (m) [matrix of math nodes,row sep=3em,column sep=5 em,minimum width=3 em] {
			\pi^*\mathbb{T}_{T^*X}	& \mathbb{T}_{\mathcal{L}^{\times}}  & pr_2^*\mathbb{T}_{\mathbb{G}_m}\\
			\pi^*\mathbb{T}_{T^*X}[-1]	& \mathbb{T}_{\mathcal{L}^{\times}}  & \mathbb{T}_{\pi}[1].  \\			
		};
		\path[-stealth]
		
		%horizontal 1st row
		(m-1-1) edge  node [above] { $ $} (m-1-2)
		(m-1-2) edge  node [above] { $ $} (m-1-3)
		%(m-1-4) edge  node [right] { } (m-1-5)
		%(m-1-5) edge  node [right] { } (m-1-6)
		%(m-1-6) edge  node [right] { } (m-1-7)
		%(m-1-7) edge  node [right] { } (m-1-8)
		%horizontal 2nd row
		(m-2-1) edge  node [above] { $  $ } (m-2-2)
		(m-2-2) edge  node [above] { $   $ } (m-2-3)
		%(m-2-4) edge  node [right] { } (m-2-5)
		%(m-2-5) edge  node [right] { } (m-2-6)
		%(m-2-6) edge  node [right] { } (m-2-7)
		%	(m-2-7) edge  node [right] { } (m-2-8)
		%vertical
		(m-1-1) edge   node [right] {$ \simeq $ } (m-2-1)
		(m-1-2) edge  node [right] { $ id $ } (m-2-2)
		(m-1-3) edge  [dashed,->] node [right] {$ \simeq $ } (m-2-3)
		%(m-1-6) edge  node [right] { } (m-2-6)
		%(m-1-7) edge  node [right] { } (m-2-7);
		%edge [dashed,-] (m-2-1)
		;
	\end{tikzpicture} 
\end{equation} Thus, we obtain an equivalence $ pr_2^*\mathbb{T}_{\mathbb{G}_m} \simeq Cone(i),$ which gives the second equivalence in (\ref{lemma_eq1}). Using the last identification, we then obtain another equivalence of  triangles

 \begin{equation}
 	\begin{tikzpicture}
 		\matrix (m) [matrix of math nodes,row sep=3em,column sep=5 em,minimum width=3 em] {
 			Cocone(\alpha) & \mathbb{T}_{\mathcal{L}^{\times}}   & Cone(i)\\
 		\pi^*\mathbb{T}_{T^*X}	& \mathbb{T}_{\mathcal{L}^{\times}}  & pr_2^*\mathbb{T}_{\mathbb{G}_m},  \\			
 			};
 		\path[-stealth]
 		
 		%horizontal 1st row
 		(m-1-1) edge  node [above] { $ i$} (m-1-2)
 		(m-1-2) edge  node [above] { $ $} (m-1-3)
 		%(m-1-4) edge  node [right] { } (m-1-5)
 		%(m-1-5) edge  node [right] { } (m-1-6)
 		%(m-1-6) edge  node [right] { } (m-1-7)
 		%(m-1-7) edge  node [right] { } (m-1-8)
 		%horizontal 2nd row
 		(m-2-1) edge  node [above] { $  $ } (m-2-2)
 		(m-2-2) edge  node [above] { $   $ } (m-2-3)
 		%(m-2-4) edge  node [right] { } (m-2-5)
 		%(m-2-5) edge  node [right] { } (m-2-6)
 		%(m-2-6) edge  node [right] { } (m-2-7)
 		%	(m-2-7) edge  node [right] { } (m-2-8)
 		%vertical
 		(m-1-1) edge  [dashed,->] node [right] {$ \simeq $ } (m-2-1)
 		(m-1-2) edge  node [right] { $ id $ } (m-2-2)
 		(m-1-3) edge  node [right] {$ \simeq $ } (m-2-3)
 		%(m-1-6) edge  node [right] { } (m-2-6)
 		%(m-1-7) edge  node [right] { } (m-2-7);
 		%edge [dashed,-] (m-2-1)
 		;
 	\end{tikzpicture} 
 \end{equation} which gives (\ref{lemma_eq2}). 
Finally, combining (\ref{lemma_eq1}) with (\ref{lemma_eq2}), we then get the homotopy cofiber sequence $ Cocone(\alpha) \rightarrow  \mathbb{T}_{\mathcal{L}^{\times}} \rightarrow Ver_{\mathcal{L}^{\times}}$ and the desired splitting in (\ref{lemma_eq3}).

\epf
%\newpage 
Now, we are in place of proving the desired result.

\begin{theorem} \label{theorem_induced contact strc1}
	Let $\mathcal{L}^{\times} \text{ and } \alpha $ be as above. Then the pair $\big(Cocone(\alpha)\hookrightarrow \mathbb{T}_{\mathcal{L}^{\times}}; \alpha\big)$ defines a 0-shifted contact structure on the derived stack $ \mathcal{L}^{\times}.$ 
\end{theorem}

\pf By Lemma \ref{lemma_useful equivalences}, we have the cofiber sequence
\begin{equation}\label{cofiber seq_horizontal_vertival}
	Cocone(\alpha)\hookrightarrow \mathbb{T}_{\mathcal{L}^{\times}} \rightarrow Ver_{\mathcal{L}^{\times}}, 
\end{equation} where $ Ver_{\mathcal{L}^{\times}}\simeq pr_2^*\mathbb{T}_{\mathbb{G}_m}.$ Since $ \mathbb{T}_{\mathbb{G}_m} \in QCoh(\mathbb{G}_m) $ is  free of rank 1, %since $\mathcal{O}(\mathbb{G}_m)=\mathbb{K}[t,t^{-1}]$ and $\Omega^1_{\mathbb{G}_m}$ is generated by $\mathrm{d}t$ over $ \mathbb{K}[t,t^{-1}]. $ Therefore, 
$ pr_2^*\mathbb{T}_{\mathbb{G}_m}$ corresponds to a line bundle $L[0]$, and hence $ \mathrm{cofib} (Cocone(\alpha)\hookrightarrow \mathbb{T}_{\mathcal{L}^{\times}}) $ is a line bundle $L[0]$ up to quasi-isomorphism. Thus, the cofiber sequence (\ref{cofiber seq_horizontal_vertival}) defines a 0-shifted pre-contact structure determined by a global 1-form $ \alpha $, with a cofiber $ Ver_{\mathcal{L}^{\times}} $ as  a certain line bundle.

Now, it remains to promote the pre-contact structure above to a 0-shifted contact structure. To this end, it suffices to show $d_{dR}\alpha$ is (locally) non-degenerate on $Cocone(\alpha)$. By the definition of $\alpha$, we get \begin{equation}
d_{dR}\alpha= \pi^*d_{dR}\lambda=\pi^*\omega,
\end{equation}which is non-degenerate on $ \pi^*\mathbb{T}_{T^*X} $ as $\omega$ is a 0-shifted symplectic structure.
Note that, from Lemma \ref{lemma_useful equivalences}, we obtain $ Cocone(\alpha)\simeq  \pi^*\mathbb{T}_{T^*X}$, which gives the desired non-degeneracy condition for $d_{dR}\alpha$ on $ Cocone(\alpha) $ and completes the proof.
\epf

\subsubsection{Prequantization of twisted cotangent stacks}\label{sec_example prequantum of twisted cotangent stacks}

Recall that  if $ X $ is a derived Artin stack locally of finite presentation equipped with a closed 1-form $\beta $ of degree $ (n + 1) $, we define the \emph{\textbf{$ n $-shifted $\beta$-twisted cotangent stack of $ X $}} to be the fiber product
\begin{equation}
	\begin{tikzpicture}
		\matrix (m) [matrix of math nodes,row sep=3em,column sep=4 em,minimum width=2 em] {
			T^*_{\beta} [n] X & X  \\
			         X  &  T^*[n+1]X , \\};
		\path[-stealth]
		(m-1-1) edge  node [left] { $   $} (m-2-1)
		(m-1-1) edge  node [above] { $   $} (m-1-2)
		%(m-1-1) edge  node [below] {} node [below] {{\small stacks}} (m-2-2)
		%(m-1-1) edge  node [below] {} node [below] {{\small  higher stacks}} (m-3-2)
		(m-2-1) edge  node  [above] {$ \Gamma_{\beta} $ } (m-2-2)
		(m-1-2) edge  node [right] {$ \Gamma_0 $ } (m-2-2);
		%edge [dashed,-] (m-2-1);
	\end{tikzpicture}
\end{equation}where $\Gamma_{\beta},\Gamma_{0}: X \rightarrow T^*[n+1]X$ are the graphs of $\beta$ and $0$-section, respectively. Note that from \cite[Corollary 2.4]{Calaque2016}, both morphisms $ \Gamma_{\beta},\Gamma_{0} $ have $(n+1)$-shifted Lagrangian structures, and hence the resulting fiber product $ T^*_{\beta} [n] X $ is an $n$-shifted symplectic stack. Moreover, \cite[Prop. 1.21]{Safronov2023} shows that the projection $T^*_{\beta} [n] X \rightarrow X$ carries a natural structure of an $n$-shifted Lagrangian fibration.

Regarding a possible prequantum structure on  $T^*_{\beta} [n] X \rightarrow X$, it has been proven in \cite[Theorem 2.24]{Safronov2023} that if $ X $ is a derived Artin stack locally of finite presentation with an $n$-gerbe $\mathcal{G}$ on $X$ such that $\beta:= c_1(\mathcal{G}) \in \mathcal{A}^{1,cl}(X,n+1)$, the \emph{characteristic class of} $\mathcal{G}$, then the projection \begin{equation*}
\pi: T^*_{c_1(\mathcal{G})} [n] X \rightarrow X
\end{equation*} has a natural structure of a prequantum $n$-shifted Lagrangian fibration determined by $\mathcal{G}$. In brief, the structure is given as follows:
\begin{observation} \label{observ_nullhomotopy of c_1 as a form in deg k-1}
Recall from \cite[Prop. 2.21]{Safronov2023} that $\pi_X: T^*[n+1] X\rightarrow X$ has a natural structure of a prequantum $ (n+1)$-shifted Lagrangian fibration, determined by the trivial $ (n+1) $-gerbe $\mathcal{G}_X$ on $X$ together with a
connective structure on $\pi_X^*\mathcal{G}_X$ given by the Liouville form $\lambda_X $. In fact, $\lambda_X \in \mathcal{A}^{1}(T^*[n+1]X,n+1)$ represents the connective structure on $\pi_X^*\mathcal{G}_X$, namely a null-homotopy of $c_1(\pi_X^*\mathcal{G}_X) \in \mathcal{A}^{1}(T^*[n+1]X,n+2)$.
Note that, from \cite[Theorem 2.9]{PTVV},  there is an induced map $$\mathcal{A}^{1,cl}(T^*[n+1]X,n+2)\rightarrow \mathcal{A}^{1,cl}(T^*_{c_1(\mathcal{G})} [n] X,n+1).$$ Now, given an $n$-gerbe $\mathcal{G}$ on $X$, \cite[Prop. 2.21]{Safronov2023} implies that the image of $ c_1(\pi_X^*\mathcal{G}_X)  $ under this map is exactly $c_1(\pi^*\mathcal{G})$. Thus, the connective structure on $ \pi_X^*\mathcal{G}_X $, given by $\lambda_X$, determines a suitable connective structure on $ \pi^*\mathcal{G}$ satisfying prequantization conditions for $\pi$. 

Let us denote the \emph{\textbf{induced connective structure}} on $ \pi^*\mathcal{G}$ by $\lambda_{c_1(\mathcal{G})}$ in $\mathcal{A}^{1}(T^*_{c_1(\mathcal{G})} [n] X,n). $
\end{observation}
\begin{example} \emph{(Construction of a shifted contact structure using a twisted cotangent stack.)}
Let $ X $ be a derived Artin stack locally of finite presentation with a trivial $0$-gerbe $\mathcal{G}$ on $X$ and a characteristic class $c_1(\mathcal{G}) \in \mathcal{A}^{1,cl}(X,1)$. 

Consider the twisted 0-shifted cotangent stack
	 $ \pi_{c_1(\mathcal{G})}: T^*_{c_1(\mathcal{G})}  X \rightarrow X $
equipped with a prequantum $0$-shifted Lagrangian fibration structure determined by $\mathcal{G}$ and a  connective structure $\lambda_{c_1(\mathcal{G})} \in \mathcal{A}^{1}(T^*_{c_1(\mathcal{G})}  X,0)$ on $ \pi^*\mathcal{G}$ described in the previous observation.

Since the trivial 0-gerbe $\mathcal{G}$ corresponds to a trivial line bundle, using the same approach in $\S$\ref{example_contact induced by prequantum on the cotangent} with obvious modifications, we define a trivial $\mathbb{G}_m$-bundle $\pi: \mathcal{L}^{\times}_{c_1(\mathcal{G})} \rightarrow T^*_{c_1(\mathcal{G})}  X$ and a global 1-form on $ \mathcal{L}^{\times}_{c_1(\mathcal{G})} $\begin{equation}\label{defn_global alpha for c_1(G)}
\alpha_{c_1(\mathcal{G})}:= \pi^* \lambda_{c_1(\mathcal{G})} + d_{dR} \log (t).
\end{equation} Then we have: 
\begin{corollary}\label{corollary_contact structure via a glabal 1-form with c_1(g)}
	Let $\mathcal{L}^{\times}_{c_1(\mathcal{G})} \text{ and } \alpha_{c_1(\mathcal{G})} $ be as above. Then the pair $\Big(Cocone(\alpha_{c_1(\mathcal{G})})\hookrightarrow \mathbb{T}_{\mathcal{L}^{\times}_{c_1(\mathcal{G})}}; \alpha_{c_1(\mathcal{G})}\Big)$ defines a 0-shifted contact structure on the derived stack $ \mathcal{L}^{\times}_{c_1(\mathcal{G})}.$ 
\end{corollary} 
\pf
The claim follows from Lemma \ref{lemma_useful equivalences} and Theorem \ref{theorem_induced contact strc1} with obvious modifications according to Observation \ref{observ_nullhomotopy of c_1 as a form in deg k-1}.
\epf
\end{example}

\subsubsection{Prequantization of the moduli stack of flat $G$-connections} \label{sec_example stack of flat conn}

Let us recall some terminology and key results from \cite[$\S$4.5]{Safronov2023}.

 Denote by $ BG  $ the \emph{classifying stack} of an affine algebraic
group $ G $ equipped with nondegenerate invariant symmetric bilinear pairing $ \langle -, - \rangle $ on its Lie algebra. More precisely, it is defined as the quotient stack
\begin{equation}
	BG= \mathrm{colim} \bigg(*\mathrel{\substack{\textstyle\leftarrow\\[-0.1ex]
			\textstyle\leftarrow \\[-0.1ex]}}  G
	\mathrel{\substack{\textstyle\leftarrow\\[-0.1ex]
			\textstyle\leftarrow \\[-0.1ex]
			\textstyle\leftarrow\\[-0.3ex]}} G \times G
	\mathrel{\substack{\textstyle\leftarrow\\[-0.1ex]
			\textstyle\leftarrow \\[-0.1ex]
			\textstyle\leftarrow \\[-0.1ex]
			\textstyle\leftarrow \\[-0.3ex]}} 
	\cdot\cdot\cdot \bigg),
\end{equation}  where the maps are given by the action and projection. Note that $ BG $ carries a canonical 2-shifted symplectic structure $\omega$. 

Fix the pair $(BG, \omega)$. Given  a smooth and proper curve $C$, we let
\begin{equation*}
	LocSys_G(C) := Map(C_{dR},BG)\footnote{Here $C_{dR}$ denotes the \emph{de Rham stack} associate with $C$. In general, for a derived stack $X$, its de Rham stack is defined to be the functor $X_{dR}: A\in cdga^{<0} \mapsto X_{dR}(A):= X(\pi_0(A)_{red})$, which corresponds to $A$-reduced points of $X$.}
\end{equation*}
be the \emph{moduli stack of flat $ G $-connections on $ C $} and
\begin{equation*}
	Bun_G(C) := Map(C,BG)
\end{equation*}
be the \emph{moduli stack of $ G$-bundles on $ C $.} Since $BG$ is 2-symplectic, PTVV's results for mapping stacks in \cite{PTVV} imply two important consequences:
\begin{enumerate}
	\item $ LocSys_G(C)$ is 0-symplectic.
	\item There is a natural closed 1-form of degree 1 on $ Bun_G(C),$ which can be obtained by \emph{integration along $C$}. We denote this form by $\int_C ev^* \omega$. For more details see \cite[$\S$1.5]{Safronov2023}.
\end{enumerate}
Note also that if $ \langle -, - \rangle $ is the Killing form, it follows from the
Grothendieck–Riemann–Roch theorem that the closed 1-form of degree 1 $\int_C ev^* \omega$ coincides with the first Chern class of the \emph{determinant line bundle} $ \mathcal{G}$ on $ Bun_G(C)$, see \cite[Example 1.26]{Safronov2023}. 

Regarding  prequantization,  we  have the identification \cite[Prop. 1.24]{Safronov2023} \begin{equation}
	LocSys_G(C)\simeq T^*_{\int_C ev^* \omega}Bun_G(C),
\end{equation}where $\int_C ev^* \omega={c_1(\mathcal{G})}.$ That is, $ LocSys_G(C) $ can be equivalently seen as a twisted cotangent stack of $ Bun_G(C)$ with the twisting 1-form $ \int_C ev^* \omega \in \mathcal{A}^{1,cl}(Bun_G(C), 1).$ From \cite[Prop. 4.22]{Safronov2023}, there is a natural prequantum 0-shifted Lagrangian fibration structure on $$  \pi_{c_1(\mathcal{G})}: LocSys_G(C) \rightarrow Bun_G(C) $$ determined by the determinant line bundle $ \mathcal{G} $ on $ Bun_G(C),$ such that $ {c_1(\mathcal{G})}= \int_C ev^* \omega,$  with a  connective structure $\lambda_{c_1(\mathcal{G})} \in \mathcal{A}^{1}(LocSys_G(C),0)$ on $ \pi^*_{c_1(\mathcal{G})}\mathcal{G}$ as described in  Observation \ref{observ_nullhomotopy of c_1 as a form in deg k-1}.

Now, we consider the $\mathbb{G}_m$-bundle $\pi: \mathcal{L}^{\times}_{\int_C ev^* \omega} \rightarrow LocSys_G(C)$ associated with $ \pi^*_{c_1(\mathcal{G})}\mathcal{G}$. If, in addition, $ \mathcal{L}^{\times}_{\int_C ev^* \omega} $ is trivial, we can define a global 1-form $ \alpha_{\int_C ev^* \omega} $ on $ \mathcal{L}^{\times}_{\int_C ev^* \omega}$  as in (\ref{defn_global alpha for c_1(G)}), which induces the desired contact structure on $ \mathcal{L}^{\times}_{\int_C ev^* \omega}$ by Corollary \ref{corollary_contact structure via a glabal 1-form with c_1(g)}. 
%\newpage

In other words, we prove:

\begin{corollary}
Let $LocSys_G(C), Bun_G(C), \mathcal{G}, \text{ and } \mathcal{L}^{\times}_{\int_C ev^* \omega} $ be as above. 
If, in addition, $\mathcal{G} $ is trivial, then the pair $$\Big(Cocone(\alpha_{\int_C ev^* \omega})\hookrightarrow \mathbb{T}_{\mathcal{L}^{\times}_{\int_C ev^* \omega}}; \ \alpha_{\int_C ev^* \omega}\Big)$$ defines a 0-shifted contact structure on the derived stack $ \mathcal{L}^{\times}_{\int_C ev^* \omega},$ where $\int_C ev^* \omega={c_1(\mathcal{G})}.$
\end{corollary}

\appendix 

\section{Darboux forms with even shifts} \label{app. non-odd shifts} For the sake of completeness, we briefly summarize the cases when $k/2$ is even or odd integer. Here, the main difference from the case $k$ being odd is about the existence of\textit{ middle degree variables}. In fact, when $k$ is odd ($k/2 \notin \mathbb{Z}$), there is no such degree. But if $k/2$ is even, then 2-forms in the middle degree variables are \textit{anti-symmetric}. On the other hand, when $k/2$ is odd, such forms are \textit{symmetric} in the middle degree variables. We directly follow \cite[Examples 5.9 \& 5.10]{Brav}.

\begin{enumerate}
	\item[(a)]  When $k/2$ is \emph{even}, say $k=-4\ell$ for $\ell\in \N$, the cdga $A$ is now free over $A(0)$ generated by the new set of variables 
	\begin{align} \label{new local variables for k=-4l}
		& x_1^{-i}, x_2^{-i}, \dots, x_{m_i}^{-i}  \text{ in degree } -i \text{ for } i= 1, 2, \dots, 2\ell-1, \nonumber \\
		& x_1^{-2\ell}, x_2^{-2\ell}, \dots, x_{m_{2\ell}}^{-2\ell}, y_1^{-2\ell}, y_2^{-2\ell}, \dots, y_{m_{2\ell}}^{-2\ell}  \text{ in degree } -2\ell, \nonumber\\
		& y_1^{k+i}, y_2^{k+i}, \dots, y_{m_i}^{k+i} \text{ in degree } k+i \text{ for } i=0,1,\dots, 2\ell-1.
	\end{align}
	We also define the element $\phi \in (\Omega^1_{A})^k$ as before.  Choose an element $H\in A^{k+1}$, the Hamiltonian, satisfying the analogue of \textit{classical master equation}, and define $d$ on $x_j^{-i},y_j^{k+i}$ (no distinguished generator $z^k$ contrary to the contact case) as in Eqn. (\ref{defn_internal d contact}) using $H$. Then $d_{dR}\alpha_0$ defines an element $
	\omega^0= \sum_{i=0}^{2\ell} \sum_{j=1}^{m_i} d_{dR}x_j^{-i} d_{dR}y_j^{k+i}$ in $ (\Lambda^2\Omega^1_{A})^k$, and set $ \omega:= (\omega^0,0,0, \cdots) $ as before, which satisfies the requirements by \cite[Example 5.9]{Brav}.
	\item[(b)]  When $k/2$ is \emph{odd}, say $k=-4\ell-2$ for $\ell\in \N$, $A$ is freely generated over $A(0)$ by the variables
	\begin{align} \label{new local variables for k=-4l-2}
		& x_1^{-i}, x_2^{-i}, \dots, x_{m_i}^{-i}  \text{ in degree } -i \text{ for } i= 1, 2, \dots, 2\ell, \nonumber \\
		& z_1^{-2\ell-1}, z_2^{-2\ell-1}, \dots, z_{m_{2\ell+1}}^{-2\ell-1} \text{ in degree } -2\ell-1, \nonumber\\
		& y_1^{k+i}, y_2^{k+i}, \dots, y_{m_i}^{k+i} \text{ in degree } k+i \text{ for } i=0,1,\dots, 2\ell.
	\end{align}
	Choose an element $H\in A^{k+1}$, the Hamiltonian, satisfying the analogue of \textit{classical master equation}
	\begin{equation} 
		\displaystyle \sum_{i=1}^{2\ell} \sum_{j=1}^{m_i} \dfrac{\partial H}{\partial x_j^{-i}} \dfrac{\partial H}{\partial y_j^{k+i}} + \frac{1}{4} \sum_{j=1}^{m_{2\ell+1}} \Big(\dfrac{\partial H}{\partial z_j^{-2\ell-1}}\Big)^2=0 \text{ in } A^{k+2}.
	\end{equation} Define $d$ on $x_j^{-i},y_j^{k+i}$ as in Equation( \ref{defn_internal d contact}) using $H$, and set $ dz_j^{-2\ell-1}:=\dfrac{1}{2}\dfrac{\partial H}{\partial z_j^{-2\ell-1}}. $ 
	
	We define the element $\phi \in \Omega^1_{A})^k$ by  \begin{equation} \label{defn_phiv2}
		\phi:=  \sum_{i=0}^{2\ell} \sum_{j=1}^{m_i} \big[-ix_j^{-i}d_{dR}y_j^{k+i} +(-1)^{i+1}(k+i)y_j^{k+i}d_{dR}x_j^i\big]+k  \sum_{j=1}^{m_{2\ell+1}} z_j^{-2\ell-1} d_{dR}z_j^{-2\ell-1}.
	\end{equation} Then $d_{dR}\alpha_0$ defines an element $
	\omega^0= \sum_{i=0}^{2\ell} \sum_{j=1}^{m_i} d_{dR}x_j^{-i} d_{dR}y_j^{k+i}+  \sum_{j=1}^{m_{2\ell+1}} d_{dR}z_j^{-2\ell-1} d_{dR}z_j^{-2\ell-1}$ in $ (\Lambda^2\Omega^1_{A})^k$, and set $ \omega:= (\omega^0,0,0, \cdots) $ as before, which satisfies the desired properties by \cite[Example 5.10]{Brav}.  
\end{enumerate}
\begin{observation}
For $k \not\equiv 2 \mod 4 $, the virtual dimension $\vdim A$ is always \emph{even}.  Otherwise, it can take any value in $ \mathbb{Z} $. More precisely,  for any $k<0$ we have 
	\[\vdim A = 
	\begin{cases*} 
		0, & \text{ if }$ k $ \mbox{is odd},\\
		\text{even in } \mathbb{Z}, & \text{ if }$ k/2 $ \mbox{is even}, \\
		\text{any value in } \mathbb{Z}, & \text{ if }$ k/2 $ \mbox{is odd}. \\
	\end{cases*}
	\]
\end{observation}

\section{Symplectic Darboux models for derived Artin stacks} Now, we give the prototype construction from \cite[Theorem 2.10]{BenBassat}:
\begin{example}  \label{example_Artin}
	Let $(\bf X, \omega)$ be a $k$-shifted symplectic derived Artin $\K$-stack and $x\in \bf X$. We construct a local model with an atlas for the case $ k $ \emph{odd}, say $k=-2\ell-1$ for $\ell\in \N.$ 
	
	We will essentially obtain either analogous or identical equations as in the case of shifted symplectic derived schemes (cf. $\S$ \ref{the pair}), but with additional finitely many generators in degree $k-1$. It means that our model will still rely on the inductively constructed sequence of cdgas as in Equation (\ref{A(n) construction}) with $A=A(-k+1)$ instead of $A(-k)$.
	
	The key idea is that the extra generators in degree $k-1$ would not play an essential role. This is because the main ingredients of the construction in  $\S$ \ref{the pair}, namely $\omega^0, H, \phi$, do not involve any of these extra variables due to degree reasons. 
	
	Applying Theorem \ref{thm_localmodelforArtincase} to $(\bf X, \omega)$, let us start with the construction of a minimal standard form open neighborhood $(A, \varphi, p)$ of $x$: Let $A(0)$ be a smooth $\mathbb{K}$-algebra of $\dim m_0$, choose $x_1^0, \dots, x_{m_0}^0$ such that $d_{dR}x_1^0, \dots, d_{dR}x_{m_0}^0$ form a basis for $\Omega_{A(0)}^1$. Then we define  $A$, as a commutative graded algebra, to be the free graded algebra over $A(0)$ generated by the variables  
	\begin{align} 
		& x_1^{-i}, x_2^{-i}, \dots, x_{m_i}^{-i}  \text{ in degree } (-i) \text{ for } i= 1, 2, \dots, \ell, \nonumber \\
		&y_1^{k+i}, y_2^{k+i}, \dots, y_{m_i}^{k+i} \text{ in degree } (k+i) \text{ for } i=0,1,\cdots \ell, \nonumber \\
		&w_1^{k-1}, w_2^{k-1}, \dots, w_{n}^{k-1} \text{ in degree } (k-1),
	\end{align} where $ m_1, \dots, m_{\ell} \in \N$ and $n=\dim H^1(\mathbb{L}_{{\bf X}}|_x)$, the (minimal possible) relative dimension of $\varphi$. Then we define an element $
	\omega^0= \sum_{i=0}^{\ell} \sum_{j=1}^{m_i} d_{dR}x_j^{-i} d_{dR}y_j^{k+i}$ in $ (\Lambda^2\Omega^1_{A})^k$, and  $ \omega:= (\omega^0,0,0, \dots) $ as before.
	
	Choose an element $H\in A^{k+1}$, the Hamiltonian, satisfying the   \textit{classical master equation} in  (\ref{defn_CME}).
	Then we define the internal differential on $A$ by $d=0$ on $A(0)$  and  by Equation (\ref{defn_internal d}). As discussed before, the condition on $H$ above is equivalent to saying $``dH=0"$. Note that we do not specify $dw^{k-1}_j$ for $j=1,\dots,n$, and hence $d$ is not completely determined on $A$ yet. But, by \cite[Theorem 2.10]{BenBassat}, $ w_1^{k-1}, w_2^{k-1}, \dots, w_{n}^{k-1} $ do not play any significant  role in the construction, and hence can be chosen arbitrarily. In fact, from the minimality argument, we have $dw^{k-1}_j|_p=0$ for each $j$.
	
	Now, choose $\phi:=  \sum_{i=0}^{\ell} \sum_{j=1}^{m_i} \big[-ix_j^{-i}d_{dR}y_j^{k+i} +(k+i)y_j^{k+i}d_{dR}x_j^i\big]$, then $dH=0$ in $A^{k+2}$, $d_{dR}H+d\phi=0$ in $(\Omega^1_{A})^{k+1}$, and $d_{dR}\phi=k\omega^0.$
	
	Let $B$ be the graded sub-cdga of $A$ over $A(0)$  generated by the variables  $x_j^i, y_j^i$ only, with inclusion $ \iota: B \hookrightarrow A. $  Since $\omega^0, H, \phi$ above do not involve any of $ w_j^{k-1} $ for degree reasons, $H\in B$, and $\omega^0, \phi$ are all images under $\iota$ of $\omega_B^0, \phi_B$, respectively. Then $ \omega_B:= (\omega^0_B,0,0, \dots) $ is a $k$-shifted symplectic structure on $V=\spec B$ such that the pair $(B, \omega_B)$ is in Darboux form as in Section \ref{the pair}, and $B$ is minimal at $p$. By construction, we have \begin{equation*}
		\spec B=V\xleftarrow{j:=\spec (\iota)} U=\spec A \xrightarrow{\varphi} \bf X
	\end{equation*} such that the induced morphism $\tau(U)\xrightarrow{\tau(j)} \tau (V)$ between truncations is an isomorphism (as $H^0(A)\simeq H^0(B)$), and $\varphi^*(\omega)\sim j^*(\omega_B)$ except in degree $k-1$. %in $k$-shifted closed 2-forms on $U$ 
	
\end{example}

\begin{remark}
	For the other cases $ (a) \ k\equiv0 \mod 4, \text{ and } (b) \ k\equiv2 \mod 4  $, the cdgas $A$ are the corresponding algebras generated by the variables as in Equations (\ref{new local variables for k=-4l}) and (\ref{new local variables for k=-4l-2}), respectively, with modification by adding the generators $ w_1^{k-1}, w_2^{k-1}, \dots, w_{n}^{k-1} \text{ in degree } k-1 $.   %We then modify $H, \phi, d$ accordingly.
\end{remark}

%\section*{Acknowledgments}
%I would like to thank Ali Ulaş Özgür Kişisel and Bayram Tekin for their useful comments and suggestions, which are helpful to clarify and improve ideas presented in this paper.
%\newpage

%-----------------------------------------------------
 
\bibliographystyle{ieeetr} % Choose Phys. Rev. style for bibliography
%Bibliography styles that can be used instead of prsty are abbrv, alpha, plain, acm, ieeetr, siam and unsrt.
\bibliography{refs}

\begin{thebibliography}{10}

\bibitem{kib}
K.~İlker Berktav, ``Shifted contact structures and their local theory,'' {\em
  arXiv:2209.09686}, Sept. 2022. To appear in the Ann. Fac. Sci. Toulouse,
  Math.

\bibitem{PTVV}
T.~Pantev, B.~To{\"e}n, M.~Vaqui{\'e}, and G.~Vezzosi, ``Shifted symplectic
  structures,'' {\em Publications math{\'e}matiques de l'IH{\'E}S}, vol.~117,
  pp.~271--328, 2013.

\bibitem{CPTVV}
D.~Calaque, T.~Pantev, B.~To\"{e}n, M.~Vaqui\'{e}, and G.~Vezzosi, ``Shifted
  {P}oisson structures and deformation quantization,'' {\em Journal of
  Topology}, vol.~10, no.~2, pp.~483--584, 2017.

\bibitem{Brav}
C.~Brav, V.~Bussi, and D.~Joyce, ``A {D}arboux theorem for derived schemes with
  shifted symplectic structure,'' {\em Journal of the American Mathematical
  Society}, vol.~32, no.~2, pp.~399--443, 2019.

\bibitem{JS}
D.~Joyce and P.~Safronov, ``A lagrangian neighbourhood theorem for shifted
  symplectic derived schemes,'' {\em Ann. Fac. Sci. Toulouse Math. 28 (2019),
  pp. 831-908}, June 2015.

\bibitem{BenBassat}
O.~Ben-Bassat, C.~Brav, V.~Bussi, and D.~Joyce, ``A `{D}arboux theorem' for
  shifted symplectic structures on derived {A}rtin stacks, with applications,''
  {\em Geometry \& Topology}, vol.~19, no.~3, pp.~1287--1359, 2015.

\bibitem{Maglio2024}
A.~Maglio, A.~G. Tortorella, and L.~Vitagliano, ``Shifted contact structures on
  differentiable stacks,'' {\em International Mathematics Research Notices},
  p.~rnae144, 2024.

\bibitem{Toen}
B.~To{\"e}n, ``Derived algebraic geometry,'' {\em EMS Surveys in Mathematical
  Sciences}, vol.~1, no.~2, pp.~153--240, 2014.

\bibitem{ToenHAG}
B.~To\"{e}n and G.~Vezzosi, ``Homotopical algebraic geometry. {II}. {G}eometric
  stacks and applications,'' {\em Memoirs of the American Mathematical
  Society}, vol.~193, no.~902, pp.~x+224, 2008.

\bibitem{Luriethesis}
J.~Lurie, {\em Derived algebraic geometry}.
\newblock PhD thesis, Massachusetts Institute of Technology, 2004.

\bibitem{Vezz2}
G.~Vezzosi, ``a derived stack?,'' {\em Notices of the AMS}, vol.~58, no.~7,
  2011.

\bibitem{Safronov2023}
P.~Safronov, ``Shifted geometric quantization,'' {\em Journal of Geometry and
  Physics}, vol.~194, p.~34, 2023.
\newblock Id/No 104992.

\bibitem{Calaque2016}
D.~Calaque, ``Shifted cotangent stacks are shifted symplectic,'' {\em arXiv
  e-prints}, p.~arXiv:1612.08101, Dec. 2016.

\bibitem{Weinstein2005}
A.~Weinstein and M.~Zambon, ``Variations on prequantization,'' in {\em Travaux
  math\'{e}matiques. {F}asc. {XVI}}, vol.~16 of {\em Trav. Math.},
  pp.~187--219, Univ. Luxemb., Luxembourg, 2005.

\end{thebibliography}

%\printbibliography %Prints bibliography when using biblatex package

\end{document}